\documentclass[12pt]{article}
\usepackage{amsmath,amsthm,amsfonts}
\usepackage{amssymb} 
\newtheorem{theorem}{Theorem}[section]
\newtheorem{proposition}[theorem]{Proposition}
\newtheorem{lemma}[theorem]{Lemma}

\newtheorem{corollary}[theorem]{Corollary}
\newtheorem{definition}[theorem]{Definition}
\newtheorem{remark}[theorem]{Remark}
\newtheorem{example}[theorem]{Example}

\newcommand{\Hom}{\mbox{\rm Hom}}

\newcommand{\id}{\mbox{\rm id}}
\renewcommand{\b}{\bullet}
\newcommand{\C}{{\cal C}}
\newcommand{\D}{{\cal D}}
\newcommand{\E}{{\cal E}}
\newcommand{\F}{{\cal F}}
\renewcommand{\H}{{\cal H}}
\newcommand{\I}{{\cal I}}
\newcommand{\J}{{\cal J}}
\newcommand{\T}{{\cal T}}
\renewcommand{\P}{{\cal P}}
\newcommand{\x}{{\mathbf x}}
\newcommand{\p}{{\mathbf p}}
\newcommand{\q}{{\mathbf q}}
\newcommand{\s}{{\mathbf s}}

\renewcommand{\O}{{\cal O}}
\renewcommand{\t}{{\mathbf t}}

\newlength{\labwidth}
\newcommand{\labarrow}[1]{
\settowidth{\labwidth}{$\scriptstyle \;\; #1 \;\;$}
\stackrel{#1}{\smash{\hbox to \labwidth{\rightarrowfill}} 
\vphantom{\longrightarrow}}
}

\begin{document}
\title{A solution of Deligne's Hochschild cohomology conjecture.}
\author{James E. McClure and Jeffrey H. Smith%
\thanks{The authors were partially supported by  NSF grants, and by SFB 343 at
the University of Bielefeld.}
\\Department of Mathematics, Purdue University, West Lafayette IN 47907--1395}
\date{February 1, 2001}
\maketitle

%mention changes from announcement

\begin{quotation}
\noindent
ABSTRACT: Deligne asked in 1993 whether the Hochschild cochain complex
of an associative ring has a natural action by the singular chains of
the little 2-cubes operad.  In this paper we give an affirmative
answer to this question.  We also show that the topological Hochschild
cohomology spectrum of an associative ring spectrum has an action by
an operad that is equivalent to the little 2-cubes operad.
\end{quotation}

\section{Introduction.}
Let us first recall some facts about
the Hochschild cochain complex $C^*(R)$ of an associative ring $R$.
An element of $C^p(R)$ is a map of abelian groups
\[
x:R^{\otimes p}\rightarrow R.
\]
Hochschild \cite{Hoch} observed that there is a cup product
in
$C^*(R)$;
if $x\in C^p(R)$ and $y\in C^q(R)$ then $x \smallsmile y$ is
the $(p+q)$-cochain defined by
\[
(x \smallsmile y) (r_1 \otimes \cdots \otimes r_{p+q}) 
=
x(r_1\otimes\cdots\otimes r_p)\cdot y(r_{p+1}\otimes\cdots\otimes r_{p+q})
\]
The cup product is clearly associative but not commutative.  In 1962
Gerstenhaber \cite{Gerst}  showed that it is chain-homotopy commutative;
the
chain homotopy between the multiplication and its twist is a sum of 
operations
\[
\circ_k:C^p(R)\otimes C^q(R)\rightarrow C^{p+q-1}(R)
\]
that take $x\otimes y$ to the $(p+q-1)$-cochain $x\circ_k y$ defined by
\[
(x\circ_k y)(r_1\otimes\cdots\otimes r_{p+q-1})=x(r_1\otimes\cdots\otimes
r_{k-1}\otimes y(r_k\otimes\cdots\otimes r_{k+q-1})\otimes\cdots\otimes
r_{p+q-1})
\]
Gerstenhaber also used the $\circ_k$ operations to construct a Lie bracket 
\[
[\ ,\ ]:C^p(R)\otimes C^q(R) \rightarrow C^{p+q-1} (R)
\]
The operations $\smallsmile$ and $[\ ,\ ]$ induce corresponding operations on the
Hochschild cohomology $H^* (R)$ which satisfy the relations making $H^*(R)$
a {\it Gerstenhaber algebra}.

Another example of a Gerstenhaber algebra is the homology of a 2-fold
loop space: $H_*(\Omega^2 A)$ (the fact that this is a Gerstenhaber
algebra is a consequence of the work of Fred Cohen \cite{FCohen}).  In
this case, the Gerstenhaber algebra structure in homology is a
consequence of the action of the little 2-cubes operad $\C_2$ on
$\Omega^2 A$. In 1993, Deligne \cite{Deligne} asked whether there was
a closer relation between these two examples: specifically, he asked
whether the Gerstenhaber algebra structure of $H^*(R)$ is induced by
an action on $C^*(R)$ of a chain operad quasi-isomorphic to the
singular chain operad of $\C_2$.  This is usually known as Deligne's
conjecture, although in the original letter it was expressed as a
desire or preference (``I would like the complex computing Hochschild
cohomology to be an algebra over [the singular chain operad of the
little 2-cubes] or a suitable version of it''). In this paper we give
an affirmative answer to this question, except that we have to replace
$C^*(R)$ by the normalized Hochschild cochains $\bar{C}^*(R)$ (a
cochain $x$ is ``normalized'' if $x(r_1\otimes\cdots\otimes r_n)=0$
whenever some $r_i=1$; see \cite[Section 10.3]{MacLane}), and in order
to make the signs work correctly it is necessary to work with the
desuspension of $\bar{C}^*$ (see \cite[page 279]{Gerst}).  Other
affirmative answers to Deligne's question for differential graded
algebras in characteristic $0$ have been found by Tamarkin \cite{T1}
and \cite{T2}, Kontsevich \cite{K}, and Voronov \cite{V}.  Tamarkin
also has a more recent approach which is similar to ours. The approach
in this paper answers Deligne's question for differential graded
algebras in any characteristic as well
as for associative ring spectra.

We show that the singular chain operad of $\C_2$ is quasi-isomorphic
(over the integers) to a suboperad of the endomorphism operad of the
desuspended reduced Hochschild cochain functor $\Sigma^{-1}\bar{C}^*$.
Explicitly, for each $n\geq 0$ let $\O(n)$ be the graded abelian
group of natural transformations
\[
\nu:(\Sigma^{-1}\bar{C}^*)^{\otimes n} \rightarrow \Sigma^{-1}\bar{C}^*
\]
(where the grading of $\nu$ is the amount by which it {\it lowers} degrees).
We give $\O(n)$ the differential 
\[
\partial(\nu)= \partial\circ\nu -\nu\circ\partial.
\]
Then $\O$ is an operad in the category of chain complexes and it is
the endomorphism operad of $\Sigma^{-1}\bar{C}^*$.  Next we
consider certain specific elements in $\O$. We have already mentioned
the element $\smallsmile\in \O(2)$.  Let
\[
e:{\mathbb Z}\rightarrow \bar{C}^*
\]
denote the element of $\O(0)$ which takes $1$ to the unit element of
the ring $R$.  For each $n\geq 2$ there is a {\it brace operation}
\[
(\bar{C}^*)^{\otimes n} \rightarrow \bar{C}^*
\]
which takes $x_1\otimes\cdots\otimes x_n$ to the cochain
\[
x_1\{x_2,\ldots,x_n\}
\]
defined by 
\[
x_1\{x_2,\ldots,x_n\}=\sum (-1)^{\varepsilon}x_1(\id,\ldots,\id,x_2,\id,\ldots,
\id,x_n,\id,\ldots,\id),
\]
where id is the identity map of $R$, 
the summation runs over all possible substitutions of $x_2$,\dots,$x_n$
(in that order) into $x_1$ and $\varepsilon$ is
\[
\sum_{j=2}^n |x_j|i_j,
\]
(here $|\ |$ denotes the degree in the desuspension and $i_j$ is the total 
number of inputs in front of $x_j$).
Note that if deg($x$)$<n$ then the sum is empty and the brace is zero.
(The brace operation for $n=2$ was defined by Gerstenhaber in \cite{Gerst}, and
the higher braces were defined by Kadeishvili \cite{Kad} and Getzler
\cite{Get}.)

Let $\H$ be the suboperad of $\O$ generated by $e$,$\smallsmile$, and the
brace operations $\{\,\}_n\in\O(n)$ for $n\geq 2$.  Our first main result is

\begin{theorem}
\label{wed1}
The singular chain operad of $\C_2$ is quasi-isomorphic as a chain operad to
$\H$.
\end{theorem}

We can give a more explicit description of $\H$ by considering the relations
satisfied by $e$, $\smallsmile$ and the $\{\,\}_n$ (\cite[Section 
1]{GV}). First we have
\begin{equation} \label{tu3}
e\smallsmile x=x\smallsmile e=x
\end{equation}
and 
\begin{equation} \label{tu4}
(x\smallsmile y)\smallsmile z=
x\smallsmile (y\smallsmile z).
\end{equation}
Because we are using normalized cochains, we have
\begin{equation} \label{tu5}
x_1\{x_2,\ldots,x_n\}=0 \quad \mbox{whenever some $x_i=1$}.
\end{equation}
There is a relation between $\{\}$ and $\smallsmile$:
\begin{equation} \label{tu6}
(x_1\cdot x_2)\{y_1,\ldots,y_n\}= 
\sum_{k=0}^n (-1)^\varepsilon x_1\{y_1,\ldots,y_k\}\cdot
x_2\{y_{k+1},\ldots,y_n\},
\end{equation}
where this time $\varepsilon=|x_2|\sum_{j=1}^k 
|y_j|$.
We can write this relation more compactly by using the standard notation 
$\nu\circ_k\nu'$ for the composition in
the operad $\O$ which inserts $\nu'$ in the $k$-th input position of $\nu$;
then equation (\ref{tu6}) becomes
\begin{equation}\label{tu666}
\{\,\}_{n+1}\,\circ_1 \smallsmile= \sum (\{\,\}_{k+1}\smallsmile\{\,\}_{n-k})\circ
\tau
\end{equation}
(where $\tau$ is the permutation that shuffles the second input to the $k+2$
position).
If we compose brace operations, we have
\begin{eqnarray} \label{tu7}
\quad x\{x_1,\ldots,x_m\}\{y_1,\ldots,y_n\}
=&\sum&
\Bigl[(-1)^\varepsilon
x\{y_1,\ldots,y_{i_1},x_1\{y_{i_1+1},\ldots,y_{j_1}\}, y_{j_1+1}, \\
&& \ldots,
y_{i_m}, x_m\{y_{i_m+1},\ldots,y_{j_m}\},y_{j_m+1},\ldots,y_m\}\Bigr],\nonumber
\end{eqnarray}
where the sum is taken over all sequences 
$0\leq i_1\leq j_1\leq i_2 \cdots \leq i_m \leq j_m\leq n$
and $\varepsilon=\sum_{p=1}^m |x_p|\sum_{q=1}^{i_p} 
|y_q|$;  this may be written more compactly as
\begin{equation}\label{tu777}
\{\,\}_{n+1} \circ_1 \{\,\}_{m+1} =
\sum\,
\{\,\}\,(\id,\ldots,\id,\{\,\},\id,\ldots,\id,\{\,\},\id,\ldots,\id)
\end{equation}
where the sum is over all ways of interpreting the right-hand side.
Since $\cup$ is a chain map we have
\begin{equation} \label{tu8}
\partial(\smallsmile)=0
\end{equation}
and the differential on a brace operation is given by
\begin{equation} \label{tu9}
\partial\{\,\}_n=-(\smallsmile\,\circ_2\,\{\,\}_{n-1})\circ\tau
+\left(\sum_i\,(\{\,\}_{n-1}\,\circ_i\,\smallsmile)\right)
-(\smallsmile\,\circ_1\,\{\,\}_{n-1})
\end{equation}
where $\tau$ is the transposition that switches the first and second entries.
If we consider $\H$ as an operad in the category of graded abelian groups 
(that is, if we neglect the differential) then $\H$ is the quotient of the 
free operad generated by $e$, $\smallsmile$ and the brace operations by the
relations 
(\ref{tu3}), 
(\ref{tu4}), 
(\ref{tu5}), 
(\ref{tu666})
and (\ref{tu777}), and the differential of $\H$ is determined by  
(\ref{tu8}) and
(\ref{tu9}).

Our method for proving Theorem \ref{wed1} is to construct a topological operad
$\C$ whose structure is based on that of $\H$ and to show that the singular
chain operad of $\C$ is quasi-isomorphic to $\H$ (this is Corollary \ref{wed2})
and that $\C$ is weakly equivalent (as an operad) to $\C_2$ (this is Theorem
\ref{mon2}).

The operad $\C$ has another significant property.  If $X^\b$ is a cosimplicial
space (or spectrum) which has cup products and $\circ_k$ operations which
satisfy the same relations as those in the Hochschild cochain complex, then
$\C$ acts on Tot$(X^\b$).  We give a careful statement of this fact in Theorem
\ref{main}.  In particular, this implies that the ``topological
Hochschild cohomology'' of an $A_\infty$ ring spectrum (see Example 
\ref{wed3} for the definition) is a $C_2$-ring spectrum.

In a sequel to this paper we will construct similar models for the little
$n$-cubes operads $\C_n$ when $n>2$.

The organization of the paper is as follows. Section \ref{sec1} is a warmup 
in which we give a sufficient condition for Tot($X^\b$) to have an $A_\infty$ 
structure.  We also introduce our basic technical tool, prismatic 
subdivision.  The results in this section are closely related to those 
in section 5 of Batanin's paper \cite{Batanin1998}; our proofs are 
simpler, but less general. In section \ref{sec2} we 
give the background needed to state 
Theorem \ref{main} and some examples to which the theorem applies.  In 
sections \ref{sec3}--\ref{sec5} we give the definition of the operad $\C$ and
prove that it acts on Tot($X^\b)$ for suitable $X^\b$. We begin in sections
\ref{sec3} and \ref{sec4} with a simplified model $\C'$ which captures most but
not all of the structure of $\C_2$, and in section \ref{sec5} we complete the
definition of  $\C$. In section \ref{sec6} we prove that the singular chain
operad of $\C$ is quasi-isomorphic to $\H$.  In sections \ref{sec7} and
\ref{sec8} we prove that $\C$ is weakly equivalent to $\C_2$, using a method of
Fiedorowicz.

We would like to thank 
Lucho Avramov,
Clemens Berger, 
Michael Brinkmeier, Dan Grayson
Rainer Vogt, and
Zig Fiedorowicz
for helpful conversations, and 
Peter May,
Jim Stasheff,
and Sasha Voronov 
for their continued interest and encouragement during the long gestation of
this project.

\newpage

\section{Maps from ${\mbox{Tot}}(X^\b)\times {\mbox{Tot}}(Y^\b)$ to 
${\mbox{Tot}}(Z^\b)$.} \label{sec1}

Suppose that $X^\bullet$, $Y^\bullet$ and $Z^\b$ are cosimplicial spaces.  
In this section we show how to construct a map 
\[
{\mbox{Tot}}(X^\b) \times {\mbox{Tot}}(Y^\b) \rightarrow {\mbox{Tot}}(Z^\b). 
\]
from cosimplicial data.

Of course, the simplest way to do this is to begin with a cosimplicial map
\[
X^\b \times Y^\b \rightarrow Z^\b
\]
and then apply Tot (using the fact that Tot commutes with products) but
we will give a construction more general than this.

In order to explain how a more general construction can be useful, let
us consider the cobar construction $X^\b$ on a based space $A$: here
$X^n$ is the Cartesian product $A^{\times n}$, the coface maps insert
basepoints and the codegeneracy maps are projections.  In this
case ${\mbox{Tot}}(X^\b)$ is homeomorphic to $\Omega A$ and thus
there is a multiplication map
\[
{\mbox{\rm Tot}}(X^\b) \times {\mbox{\rm Tot}}(X^\b) \rightarrow {\mbox{\rm Tot}}(X^\b);
\] 
in fact there is one such map for each partition of the unit 
interval into two parts.  On the other hand, there is clearly no sensible 
way to map $A^{\times n} \times A^{\times n}$ to $A^{\times n}$, so there is 
in general no cosimplicial map $X^\b \times X^\b \rightarrow X^\b$ that can 
induce the multiplication.

The following definition describes the structure at the cosimplicial level 
that our construction accepts as input; it is closely related to Batanin's
definition of the tensor product of cosimplicial spaces \cite{Batanin1993}.

\begin{definition} 
\label{sun4}
{\rm {\bf (i)} A {\it cup-pairing} $\phi:(X^\b,Y^\b)\rightarrow Z^\b$ is a 
family of maps
\[
\phi_{p,q}: X^p\times Y^q \rightarrow Z^{p+q}
\]
satisfying
\begin{description}
\item[\rm(a)]
$d^i\phi_{p,q}(x, y)=\left\{
\begin{array}{ll}
\phi_{p+1,q}(d^i x , y) & \mbox{if $i\leq p$} \\
\phi_{p,q+1}(x, d^{i-p}y) & \mbox{if $i>p$}
\end{array}
\right.
$
\item[\rm(b)]
$\phi_{p+1,q}(d^{p+1}x, y)=\phi_{p,q+1}(x, d^0 y)
$
\item[\rm(c)]
$s^i\phi_{p,q}(x, y)=\left\{
\begin{array}{ll}
\phi_{p-1,q}(s^i x, y) & \mbox{if $i\leq p-1$} \\
\phi_{p,q-1}(x , s^{i-p}y) & \mbox{if $i \geq p$}
\end{array}
\right.
$
\end{description}

{\bf (ii)} A {\it morphism} of cup pairings, from 
$\phi:(X^\bullet, Y^\bullet) \rightarrow Z^\bullet$ to
$\phi':({X'}^\bullet, {Y'}^\bullet) \rightarrow {Z'}^\bullet$ 
is a triple of cosimplicial maps
\[
\mu_1:X\rightarrow X',
\mu_2:Y\rightarrow Y',
\mu_3:Z\rightarrow Z'
\]
such that
\[
\mu_3\circ\phi'_{p,q}=\phi_{p,q}\circ (\mu_1 \times \mu_2)
\]
for all $p,q$.

{\bf (iii)} Given a cup-pairing $\phi:(X^\b,X^\b)\rightarrow X^\b$, a {\it 
unit} for
$\phi$ is a sequence of points $e_n\in X^n$ such that the set $\{\, e_n \,\}$
is closed under all cofaces and codegeneracies and 
$\phi_{0,p}(e_0,x)=\phi_{p,0}(x,e_0)=x$
for all $x$ (i.e., the $e_n$ determine a map from the trivial cosimplicial
space to $X^\b$, and $e_0$ is a unit in the usual sense).
}
\end{definition}

From now on we shall usually drop the subscripts and just write $\phi(x,y)$ for
$\phi_{p,q}(x,y)$.

\begin{remark} 
\label{sat3}
{\rm {\bf (i)} The cobar construction on $A$ has the cup-pairing 
\[
\phi((a_1,\ldots,a_p),(b_1,\ldots,b_q))=(a_1,\ldots,a_p,b_1,\ldots,b_q)
\]

{\bf (ii)} The definition of cup-pairing is modeled on 
the properties of the cup product in the Hochschild cohomology complex. 
It is easy to check that the function $\phi(x,y)=x\smallsmile y$ is a
cup-pairing in our sense if the Hochschild complex is given the usual
cosimplicial structure in which $s^i$ inserts a unit in the $i+1$-st position 
and
$d^i$ is defined by
\[
(d^i x)(r_1,\ldots,r_{p+1})=
\left\{
\begin{array}{ll}
r_1 x(r_2,\ldots,r_{p+1}) & \mbox{if $i=0$} \\
x(\ldots,r_{i}r_{i+1},\ldots) & \mbox{if $0<i<p+1$} \\
x(r_1,\ldots,r_p)r_{p+1} & \mbox{if $i=p+1$}
\end{array}
\right.
\]

{\bf (iii)} A cosimplicial map $\mu:X^\b\times Y^\b \rightarrow Z^\b$ induces a cup
pairing by 
\[
\phi(x , y) =\mu(d^{p+q}d^{p+q-1}\cdots d^{p+1}x,d^0d^0\cdots d^0 y)
\]
(that is, the last coface map is applied $q$ times to $x$, and the zeroth
coface map is applied $p$ times to $y$, where $p$ is the degree of $x$ and $q$
is the degree of $y$; then $\mu$ is applied to the resulting pair).

}\end{remark}

\begin{theorem} \label{sun1} {\bf (i)} A cup-pairing 
\[
\phi:X^\b \times Y^\b\rightarrow Z^\b
\]
induces a map
\[
\bar{\phi}_u:{\mbox{\rm Tot}}(X^\b) \times {\mbox{\rm Tot}}(Y^\b) \rightarrow {\mbox{\rm Tot}}(Z^\b). 
\]
for each $u$ with $0<u<1$.

{\bf (ii)}
A morphism of cup-pairings induces a commutative diagram
\[
\begin{array}{ccc}
{\mbox{\rm Tot}}(X^\b) \times {\mbox{\rm Tot}}(Y^\b) &\labarrow{\bar{\phi}_u} &  {\mbox{\rm Tot}}(Z^\b) \\
\downarrow && \downarrow\\
{\mbox{\rm Tot}}({X'}^\b) \times {\mbox{\rm Tot}}({Y'}^\b) &\labarrow{\bar{\phi'}_u} &  {\mbox{\rm Tot}}({Z'}^\b) \\
\end{array}
\]

{\bf (iii)} If the cup pairing comes from a cosimplicial map
$\mu$ then each $\bar{\phi}_u$ is homotopic to the usual map induced by $\mu$.
\end{theorem}

\noindent{\bf Remark.}
Part (i) of this theorem is implicit in Batanin's paper \cite{Batanin1998},
particularly in the proof of \cite[Theorem 5.2]{Batanin1998}, which situates
this result in a more general context.

\medskip

Our next result refers to ${\mbox{\rm Tot}}'$, which is a construction
related to Tot in the same way that the Moore loop space is related to
the ordinary loop space; we will give the definition below.

\begin{theorem}
\label{sun2}
A strictly associative cup-pairing on $X^\b$ with a unit induces a
strictly associative multiplication on ${\mbox{\rm Tot}}'(X^\b)$ and an
action of the little 1-cubes operad on ${\mbox{\rm Tot}}(X^\b)$.
\end{theorem}

\begin{remark} {\rm  
{\bf (i)} Batanin \cite[Theorems 5.1 and 5.2]{Batanin1998} uses trees to
construct an 
$A_\infty$ operad which acts on Tot($X^\b$) when $X^\b$ has a strictly 
associative cup-pairing. 
The $0$-th space of Batanin's operad is not a 
point, so that an action of his operad on a space or spectrum provides a 
multiplication with a homotopy unit rather than a strict unit.  

\medskip

{\bf (ii)} If $X^\bullet$ is the cobar 
construction on a space $A$ then the action of the little 1-cubes on ${\mbox{\rm Tot}}(X^\b)$
is homeomorphic to the usual action of the little 1-cubes on $\Omega(A)$.  Also, there is a continuous bijection from 
${\mbox{\rm Tot}}'(X^\b)$ to the Moore loop
space of $A$ which takes the multiplication on ${\mbox{\rm Tot}}'(X^\b)$ 
to the usual multiplication on the Moore loop space.
}\end{remark}

Before giving the proof of Theorems \ref{sun1} and \ref{sun2}, we will 
describe a way of 
subdividing a simplex which we call the ``prismatic
subdivision.'' This is a homeomorphism discovered independently 
and at various times by Lisica and Mardesic \cite{LM}, Grayson 
\cite{Grayson}, and ourselves.
(There is a related ``edgewise subdivision'' discovered by Quillen, Segal, 
B\"okstedt and Goodwillie; the edgewise subdivision of a simplex is a 
subdivision of the prismatic subdivision.)

It will be convenient from now on to let $\s$, $\t$, etc.\ stand 
for a point of a simplex and let $s_0,s_1,\ldots$ (respectively,
$t_0,t_1,\ldots$, etc.) be its coordinates.

Define a space $D^n$ for each $n \geq 0$ by
\[
D^n=
(\coprod_{p=0}^{n}
\Delta^p \times \Delta^{n-p} )/\sim,
\]
where $\sim$ is defined by
$(d^{p+1}\s, \t)\sim (\s, d^0 \t)$
if $\s\in 
\Delta^p $ and $\t \in \Delta^{n-p-1}$.
For each $u$ with $0\leq u\leq 1$, there is a map
%do I have to allow \leq instead of >? YES---for part (iii)
\[
\sigma^n(u):D^n
\rightarrow
\Delta^{n}
\]
whose restriction to
$\Delta^p \times \Delta^{n-p}$ takes 
\[
(s_0,\ldots,s_p),(t_0,\ldots,t_{n-p})
\]
to
\[
(us_0,\ldots,us_{p-1},us_p+(1-u)t_0, (1-u)t_1,\ldots,(1-u)t_{n-p}).
\]
It is easy to check that this map is well-defined, continuous and onto, 
and when $0<u<1$ it is also one-to-one.  Thus we have: 

\begin{proposition}
\label{sun3}
$\sigma^n(u)$ is a homeomorphism if $0<u<1$.
\end{proposition}

This homeomorphism is the prismatic subdivision of the simplex $\Delta^n$.

Here are the prismatic subdivisions of the 1-simplex and the
2-simplex (with $u=\frac{1}{2}$):

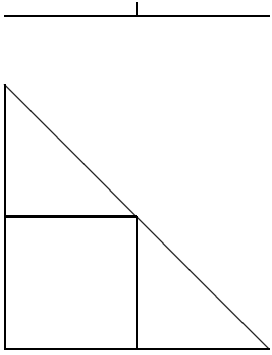
\begin{picture}(100,50)

\put(0,25){\line(1,0){100}}
\put(50,25){\line(0,1){5}}

\end{picture}

\begin{picture}(100,100)

\put(0,100){\line(1,-1){100}}
\put(0,0){\line(0,1){100}}
\put(0,0){\line(1,0){100}}
\put(0,50){\line(1,0){50}}
\put(50,0){\line(0,1){50}}

\end{picture}

\bigskip

We can now prove part (i) of Theorem \ref{sun1}
Recall that a point in ${\mbox{\rm Tot}}(X^\bullet)$ is a sequence 
\[
a_0\in X^0, a_1:\Delta^1\rightarrow X^1, a_2:\Delta^2\rightarrow X^2, \ldots
\]
which is {\it consistent}, i.e.,
\[
d^i \circ a_n = a_{n+1}\circ d^i
\]
and
\[
s^i \circ a_n = a_{n-1} \circ s^i
\]
Thus what is required is: given consistent sequences 
$a_n:\Delta^n\rightarrow X^n$
and $b_n:\Delta^n\rightarrow Y^n$, to construct a consistent sequence
$c_n:\Delta^n \rightarrow Z^n$.  

First observe that, by part (b) of Definition \ref{sun4}(i), the maps 
\[
\phi\circ(a_p \times b_{n-p}): \Delta^p \times \Delta^{n-p} 
\rightarrow Z^n
\]
fit together to give a well-defined map
\[
D^n \rightarrow Z^n
\]
(where $D^n$ is the space defined before Proposition \ref{sun3}).
We define $c_n$ for each $n$ to be the composite of this map with 
$(\sigma^n(u))^{-1}$.
The fact that the $c_n$ commute with coface and codegeneracy maps
is immediate from the formula for $\sigma^n(u)$ and parts (a) and (c) of 
Definition \ref{sun4}(i).  This concludes the proof of part (i).

Here are pictures of $c_1$ and $c_2$ (with $u=\frac{1}{2}$):

\begin{picture}(150,60)

\put(0,25){\line(1,0){150}}
\put(75,25){\line(0,1){5}}
\put(0,25){\makebox(75,20){$a_0\smallsmile b_1$}}
\put(75,25){\makebox(75,20){$a_1\smallsmile b_0$}}

\end{picture}

\begin{picture}(150,150)

\put(0,150){\line(1,-1){150}}
\put(0,0){\line(0,1){150}}
\put(0,0){\line(1,0){150}}
\put(0,75){\line(1,0){75}}
\put(75,0){\line(0,1){75}}
%\put(0,15){\makebox(60,10){$a_1\smallsmile b_1$}}
\put(0,35){\makebox(75,10){$a_1\smallsmile b_1$}}
\put(0,90){\makebox(50,10){$a_2\smallsmile b_0$}}
\put(75,15){\makebox(50,10){$a_0\smallsmile b_2$}}

\end{picture}

Part (ii) of Theorem \ref{sun1} is immediate from the definitions.

We now turn to the proof of part (iii).  By part (ii), it suffices to 
consider the case where $\mu$ is the identity map of 
$X^\b \times Y^\b$.  Let $\phi:(X^\b,Y^\b)\rightarrow X^\b\times Y^\b$ denote 
the pairing induced by the
identity map of $X^\b \times Y^\b$ (as in Remark \ref{sat3}2(iii)): thus 
\[
\phi(x,y)=
(d^{p+q}d^{p+q-1}\cdots d^{p+1}x,d^0d^0\cdots d^0 y)
\] 
Let 
\[
F: {\mbox{\rm Tot}}(X^\b) \times {\mbox{\rm Tot}}(Y^\b) \rightarrow {\mbox{\rm
Tot}}(X^\b) \times {\mbox{Tot}}(Y^\b) 
\]
denote $\bar{\phi}_{1/2}$; it suffices to show that $F$ is homotopic
to the identity.  Let $\pi_1$ and $\pi_2$ be the projections of
${\mbox{\rm Tot}}(X^\b) \times {\mbox{\rm Tot}}(Y^\b) $ on its first and
second factors; it suffices to show that $\pi_i \circ F$ is homotopic
to $\pi_i$ for $i=1, 2$. We consider the case $i=1$; the other case is
similar.  Let $\tilde{D}^n$ denote the space
\[
(\coprod_p \Delta^p \times (\Delta^{n-p} \wedge [\frac{1}{2},1]))/\sim ,
\]
where the basepoint of $[\frac{1}{2},1]$ is taken to be 1 and $\sim$ is defined
by $(d^{p+1}\s, \t \wedge u)\sim (\s, d^0 \t \wedge u)$
if $\s\in 
\Delta^p $ and $\t \in \Delta^{n-p-1}$.
The map
$\tau:\tilde{D}^n \rightarrow \Delta^n \times [\frac{1}{2},1]$
given by 
\[
\tau(\s,\t \wedge u)=(\sigma^n(u)(\s,\t),u)
\]
is a homeomorphism, since it is 1-1, onto, and continuous.  
For each $n$, the maps
\[
\Delta^p \times (\Delta^{n-p}\wedge [\frac{1}{2},1]) \times {\mbox{\rm Tot}}(X^\b) \times
{\mbox{\rm Tot}}(Y^\b) \rightarrow X^n
\]
which take
$(\s,\t\wedge u,\{a_n\},\{b_n\})$ to $a_n(d^n d^{n-1}\cdots d^{p+1}\s)$
fit together to give a well-defined map
\[
\tilde{D}^n \times {\mbox{\rm Tot}}(X^\b) \times {\mbox{\rm Tot}}(Y^\b) \rightarrow  X^n. 
\]
Composing this with $(\tau_n)^{-1}$ gives a map
\[
\Delta^n \times [\frac{1}{2},1] \times {\mbox{\rm Tot}}(X^\b) \times {\mbox{\rm
Tot}}(Y^\b) \rightarrow X^n,
\]
whose adjoint is a map 
\[
H_n: [\frac{1}{2},1] \times {\mbox{\rm Tot}}(X^\b) \times {\mbox{\rm Tot}}(Y^\b) \rightarrow 
{\mbox{M}\mbox{a}\mbox{p}}(\Delta^n,X^n).
\]
Taken together, the $H_n$ give a map
\[
H: [\frac{1}{2},1] \times {\mbox{\rm Tot}}(X^\b) \times {\mbox{\rm Tot}}(Y^\b) \rightarrow {\mbox{\rm Tot}}(X^\b)
\]
which is equal to $\pi_1 \circ F$ when $u=\frac{1}{2}$ and to $\pi_1$ when
$u=1$.  Thus $\pi_1 \circ F$ is homotopic to $\pi_1$ as required. \qed

Before proving Theorem \ref{sun2} we define ${\mbox{\rm Tot}}'$.   
First we define
$\Delta^n_p$,
for each $p\geq 0$,
to be the set
\[
\{\, (s_0,\ldots,s_n) \,|\, s_i \geq 0 \ \mbox{for all $i$ and}\
s_0+\cdots+s_n=
p \,\}
\]
(Note that $\Delta^\b_p$ is a cosimplicial space for each $p$, and that
$\Delta^\b_0$ is is the trivial cosimplicial space consisting of a point in
each degree.)
Then we define
${\mbox{\rm Tot}}'(X^\b)$ to consist of pairs $(p,\{a_n\})$, where $p$ is $\geq 0$ and
$\{a_n\}$ is a cosimplicial map $\Delta^\b_p\rightarrow X^\b$.
In order to describe the topology of ${\mbox{\rm Tot}}'$ we first observe that for $p>0$ 
the cosimplicial spaces $\Delta^\b_p$ and $\Delta^\b$ are isomorphic, and thus
there is a 
bijection between ${\mbox{\rm Tot}}'(X^\b)$ and the disjoint union
\[
\mbox{Const}(X^\b)
\cup 
({\mbox{\rm Tot}}(X) \times {\mathbb R}_{>0}),
\]
where 
$\mbox{Const}(X^\b)$ denotes 
the space of cosimplicial maps from $\Delta^\b_0$ to $X^\b$
and ${\mathbb R}_{>0}$ denotes the positive reals.  This in turn implies that
there is a bijection between $T'(X^\b)$ and the pushout
\[
\begin{array}{ccc}
\mbox{Const}(X^\b)\times {\mathbb R}_{>0} & \rightarrow & 
\mbox{Const}(X^\b)\times {\mathbb R}_{\geq 0} \\
\downarrow && \\
{\mbox{\rm Tot}}(X^\b)\times {\mathbb R}_{>0} & & 
\end{array}
\]
We give ${\mbox{\rm Tot}}'(X^\b)$ the topology which makes this bijection a homeomorphism.
(We could have simply defined ${\mbox{\rm Tot}}'(X^\b)$ to be this pushout, but this would
make the formulas in the rest of the paper more complicated).  Note that there
is a continuous projection $\rho:{\mbox{\rm Tot}}'(X^\b)\rightarrow {\mbox{\rm
Tot}}(X^\b)$ which makes
${\mbox{Tot}}(X^\b)$ a deformation retract of ${\mbox{Tot}}'(X^\b)$.

\bigskip
\noindent
{\bf Proof of Theorem \ref{sun2}.}  We need a slight generalization of the idea of 
prismatic subdivision:
define a space $D^n_{p,q}$ for each $n \geq 0$ and each $p,q\geq 0$ by
\[
D^n_{p,q}=
(\coprod_{k=0}^{n}
\Delta^k_p \times \Delta^{n-k}_q )/\sim,
\]
where $\sim$ is defined by
$(d^{k+1}\s, \t)\sim (\s, d^0 \t)$
if $\s\in 
\Delta^k_p $ and $\t \in \Delta^{n-k-1}_r$.
Define
\[
\sigma^n_{p,q}:D^n_{p,q}
\rightarrow
\Delta^{n}_{p+q}
\]
by 
\[
\Sigma((s_0,\ldots,s_k),(t_0,\ldots,t_{n-k}))
=
(s_0,\ldots,s_{k-1},s_k+t_0, t_1,\ldots,t_{n-k}).
\]
Then $\sigma^n_{p,q}$ is a homeomorphism for every choice of $p$ and $q$.

Now we can define the multiplication on ${\mbox{\rm Tot}}'$. Given a
pair of points $(p,\{a_n\})$ and $(q,\{b_n\})$ in ${\mbox{\rm Tot}}'$ then,
by part (b) of Definition \ref{sun4}(i), the maps
\[
\phi\circ(a_k \times b_{n-k}): \Delta^k_p \times \Delta^{n-k}_q 
\rightarrow X^n
\]
fit together to give a well-defined map
\[
D^n_{p,q} \rightarrow X^n
\]
Composing this map with 
$(\sigma^n_{p,q})^{-1}$ gives a map
$c_n:\Delta^n_{p+q}\rightarrow X^n$ for each $n$. 
The $c_n$ commute with coface and codegeneracy maps
by parts (a) and (c) of 
Definition \ref{sun4}(i), and we define the product of the points
$(p,\{a_n\})$ and $(q,\{b_n\})$ to be the point $(p+q,\{c_n\})$.
This multiplication is clearly associative and unital, and
we leave it to the reader to check that it is continuous.

It remains to give the action of the little 1-cubes operad $\C_1$ on 
${\mbox{\rm Tot}}(X^\b)$.  So given a point $z$ of $\C_1(k)$ and $k$ points
$a_1,\ldots,a_k$ of ${\mbox{\rm Tot}}(X^\b)$, we are required to define a point 
$\gamma(z,a_1,\ldots,a_k)$ in ${\mbox{\rm Tot}}(X^\b)$.  
The idea is to imitate the way that $\C_1$ acts on a loop space, that is, we
scale the $a$'s to the lengths of the corresponding little intervals and fill
in the blank spaces with appropriately scaled copies of the unit element.

First observe that the point $z$ of $\C_1(k)$ can
be written as a sequence of $2k+1$ numbers $(p_1,p_2,\ldots,p_{2k+1})$, where
the even numbered $p$'s are the lengths of the little intervals and the odd
numbered $p$'s are the lengths of the empty spaces.  

Next we introduce some notation.
Given a function $f:\Delta^n\rightarrow X^n$ and a number $p>0$, let us write 
$p\cdot f$ for the function 
$\Delta^n_p\rightarrow X^n$ which takes $\s$ to $f(\frac{1}{p}\cdot\s)$.  
Given a
point $a=\{a_n\}$ in ${\mbox{\rm Tot}}(X^\b)$, let us write 
$p\cdot a$ for the point $(p,\{p\cdot a_n\})$ in ${\mbox{\rm Tot}}'(X^\b)$.  Recall that the
cup product has a unit $e$, which is a map from $\Delta^\b_0$ to $X^\b$, and
for each $p \geq 0$ write $e_p$ for the composite 
\[
\Delta^\b_p \rightarrow \Delta^\b_0 \labarrow{e} X^\b
\]
Let $*$ denote the multiplication on ${\mbox{\rm Tot}}'(X^\b)$.  Then we define 
$\gamma(z,a_1,\ldots,a_k)$  to be
\[
\rho(e_{p_1} *(p_2\cdot a_1) *e_{p_3} *(p_4\cdot a_2)*\cdots* e_{p_{2k+1}}),
\]
where $\rho$ is the projection from ${\mbox{\rm Tot}}'(X^\b)$ to ${\mbox{\rm Tot}}(X^\b)$.
\qed

\newpage

\section{Operads with multiplication and their associated cosimplicial objects.}
\label{sec2}

In this section we give a formal statement of the fact that if $X^\b$ is a
cosimplicial space or spectrum with $\smallsmile$ and $\circ_k$ operations that
satisfy the same relations as those in the Hochschild cochain complex then our
operad $\C$ acts on Tot($X^\b$).  For this we need some background.

First recall that if $\O$ is an operad with structure maps
\[
\gamma:\O(n)\times\O(j_1)\times\cdots\times\O(j_n)\rightarrow 
\O(j_1+\cdots +j_n)
\]
and identity element ${\mbox{id}}\in\O(1)$, it is customary to write 
$\circ_i$ for the map
\[
\O(n)\times\O(j)\rightarrow \O(n+j-1)
\]
which takes $(o_1,o_2)$ to 
$\gamma(o_1,\mbox{id},\ldots,o_2,\ldots,\mbox{id})$
(with $i-1$ id's before the $o_2$).

Next let $R$ be a ring and let 
$\O_R$ denote the endomorphism operad of
the underlying abelian group of $R$; that is, $\O_R(n)$ is the set of 
homomorphisms of abelian groups from $R^{\otimes n}$ to $R$, with the evident
operad structure.  There are special elements
$\mu\in \O_R(2)$ (the multiplication in the ring $R$) and $e\in \O_R(0)$
(the unit element of $R$).  Of course, $\O_R(n)$ is the same as the
group of $n$-cochains in the Hochschild complex of $R$, and the cosimplicial 
structure of the Hochschild complex can be recovered from the operad structure
of $\O_R$ and the elements $\mu$ and $e$:
\begin{equation} \label{e1}
(d^i x)=
\left\{
\begin{array}{ll}
\mu\circ_2 x & \mbox{if $i=0$} \\
x\circ_i \mu & \mbox{if $0<i<p+1$} \\
\mu\circ_1 x & \mbox{if $i=p+1$}
\end{array}
\right.
\end{equation}
\begin{equation} \label{e2}
s_i(x) = x\circ_{i+1} e
\end{equation}
When the $d^i$ and $s^i$ are defined in this way the 
cosimplicial identities follow formally from the identities 
\begin{equation} \label{e3}
\mu\circ_1 \mu=\mu\circ_2 \mu 
\end{equation}
and
\begin{equation} \label{e4}
\mu\circ_1 e=\mu\circ_2 e=\mbox{id}.
\end{equation}

This example motivates the following definition, which is due to Gerstenhaber
and Voronov \cite{GV}.  Recall (\cite[Definition 3.12]{MayG}) that a
non-$\Sigma$ operad is a structure which has all the properties of an operad
except those having to do with the actions of the symmetric groups.  Also
recall that the definition of operad makes sense in any symmetric monoidal
category.

\begin{definition}
{\rm
Let ${\cal S}$ be a symmetric monoidal category with unit object $S$.
An {\it operad with multiplication} in $\cal S$
is a non-$\Sigma$ operad $\O$ in $\cal S$ together with
maps $e:S\rightarrow \O(0)$ and $\mu:S\rightarrow \O(2)$ satisfying (\ref{e3}) 
and (\ref{e4}).  
The {\it associated cosimplicial object} $\O^\b$ consists of the objects 
$\O(n)$ with coface and codegeneracy maps given by (\ref{e1}) and (\ref{e2}).
}
\end{definition}

The reader will perhaps be relieved to know that the only symmetric monoidal
categories we will be concerned with in this paper are the category of spaces,
the category of $S$-modules \cite[Chapter II]{EKMM}, and (briefly) the category
of sets.

\begin{remark}
{\rm
{\bf (i)}
If we let $Ass$ be the non-$\Sigma$ operad for which $Ass(n)$
is the object $S$ for each $n$ (so that an algebra over $Ass$ is an associative
monoid in $\cal S$) then a more economical (but equivalent) way to define an 
operad with multiplication is: a non-$\Sigma$ operad together with a 
non-$\Sigma$-operad morphism from $Ass$.

{\bf (ii)}
If $\O$ is an operad with multiplication then the associated cosimplicial
object has a cup-pairing in the sense of Section \ref{sec1}: we define 
$\phi_{p,q}(x,y)$
to be $(\mu\circ_1 x)\circ_{p+1} y$.
}
\end{remark}

We can now state our main result.  Recall that a weak equivalence of operads is
an operad morphism which is a weak equivalence on each object.  We say that two
operads are {\it weakly equivalent} if there is a third operad which maps to
each of them by a weak equivalence.

\begin{theorem} \label{main}
There is an operad $\C$ in the category of spaces with the 
following properties:

{\bf (i)} $\C$ is weakly equivalent to the little 2-cubes operad $\C_2$.

{\bf (ii)} If $\O$ is an operad with multiplication in the category of spaces
or spectra then $\C$ acts on ${\mbox{\rm Tot}}(\O^\b)$.
\end{theorem}

Here are some examples to which Theorem \ref{main} can be applied.  In the
examples that refer to the category of $S$-modules 
we will always write $\wedge$ for $\wedge_S$.

\begin{example} \label{wed3}
{\bf (The Hochschild cohomology complex of a ring spectrum.)}
{\rm
Let $R$ be an $S$-algebra in the sense of \cite[Section II.3]{EKMM} and let
$\O$ be the endomorphism operad of $R$:
\[
\O(n)=F_S(R^{\wedge n}, R)
\] 
(see \cite[Section II.1]{EKMM} for the definition of $F_S$).
Let $e:S\rightarrow R$ be the
unit map of $R$ and let $\mu:S\rightarrow \O(2)$ be adjoint to the multiplication
map $R\wedge R\rightarrow R$.  Then $\O$ is an operad with multiplication and
we define the Hochschild cohomology complex of $R$ to be the cosimplicial
$S$-module associated to $\O$. 
}
\end{example}

\begin{example}[The loop space of a topological monoid.] \label{ex2}
{\rm (We would like to thank Zig Fiedorowicz for pointing out this example to
us.) 
Let $A$ be a topological monoid, with the unit of $A$ chosen as basepoint. 
We will define an operad with multiplication
$\O$ whose associated cosimplicial space is the cobar construction on $A$.
Of course, we let $\O(n)=A^{\times n}$.  In order to define the operad 
structure we observe (as in \cite[p.\ 6]{MayOp})
%this is May's paper in Contemp Math 202
that it suffices to specify the operations 
\[
\circ_i: \O(n)\times\O(j)\rightarrow \O(n+j-1),
\]
and we define these by the equation
\begin{equation} \label{e5}
(a_1,\ldots,a_n)\circ_i (b_1,\ldots,b_j)=
(a_1,\ldots,a_{i-1},a_i b_1,\ldots,a_i b_j, a_{i+1},\ldots a_n)
\end{equation}
Finally, we choose $e\in\O(0)$ and $\mu\in \O(2)$ to be the basepoints.
(Observe that the cup-pairing
determined by $\mu$ is the same as that defined in Remark \ref{sat3}(i).)
}
\end{example}

In preparation for our next two examples it is helpful to reformulate Example 
\ref{ex2} in a fancier way.  Let $B_\b$ denote the simplicial circle
$\Delta^1/\partial\Delta^1$. 
Recall that a non-basepoint simplex in $B_n$
is a sequence $(b_0,b_1,\ldots,b_n)$ of $n+1$ zeroes and ones, with 
the zeroes coming before the ones.
We can put the non-basepoint simplices of $B_n$ in one-to-one correspondence
with the numbers 1 to $n$ by letting the number $j$ correspond to the simplex
with $j$ zeroes, and this gives a homeomorphism $\alpha$ between $\O(n)$ and 
the space of based maps from $B_n$ to $A$; under $\alpha$ the
coface and codegeneracy maps of $\O(n)$ agree with the maps induced by the
faces and degeneracies of $B_n$. If $x\in\O(n)$ and $b\in B_n$ then we write 
$x(b)$ for the element $\alpha(x)(b)$ of $A$.
Now equation (\ref{e5}) is equivalent to the following:
if $x\in \O(n)$, $y\in \O(j)$ and $b\in B_n$ then
\begin{equation} \label{e6}
(x\circ_i y)(b)=x({b'_{ij}})\cdot y({b''_{ij}}),
\end{equation}
where $\cdot$ denotes the product in $A$, $b'_{ij}$ is the simplex
\[
(b_0,b_1,\ldots,b_{i-1},b_{i+j-1},\ldots,b_{n+j-1})
\]
and $b''_{ij}$ is the simplex
\[
(b_{i-1},\ldots,b_{i+j-1}).
\]

We can generalize the definitions of $b'_{ij}$ and $b''_{ij}$ to any simplicial
set $B_\b$ as follows.
If $s\leq t$
we write $\partial(s,t)$ for the composite of face maps
\[
\partial_{s}\partial_{s+1}\cdots\partial_{t}
\]
Now if $n\geq 1$, $j\geq 0$, $1\leq i\leq n$, and $b\in B_{n+j-1}$,  we define
\[
b'_{ij}=\partial(i,i+j-2) b
\]
and 
\[
b''_{ij}=
\partial(0,i-2)\partial(i+j,n+j-1)b.
\]
These definitions will be familiar to many readers because they are related to
Steenrod's original definition of the $\smallsmile_1$ product 
\cite{Steenrod}: if $\xi\in C^n(B_\b;{\mathbb Z}/2)$
and $\eta\in C^j(B_\b;{\mathbb Z}/2)$ are mod-2 simplicial cochains of $B_\b$ and
$b\in B_{n+j-1}$ then
\[
(\xi\smallsmile_1 \eta)(b)=
\sum_i \xi(b'_{ij})\eta(b''_{ij})
\]

\begin{example}[$\Omega^k$ of a space when $k\geq 2$.] \label{ex3}
{\rm 
Let $A$ be any based space and fix $k\geq 2$. We can define a ``higher'' 
cobar construction $X^\b$, with ${\mbox{\rm Tot}}(X^\b)$ homeomorphic to $\Omega^k 
A$, as follows:
let $B_\b$ be the simplicial $k$-sphere $\Delta^k/\partial 
\Delta^k$ and let $X^n$ be the space of based maps 
from $B_n$ to $A$, with coface and codegeneracy maps induced by the face and
degeneracy maps of $B_\b$.  
Next we define an operad with multiplication $\O$ whose associated
cosimplicial space is $X^\b$.  Of course, we let $\O(n)=X^n$.  
We define the $\circ_i$ operations by equation (\ref{e6}), but note that the
symbol $\cdot$ now requires interpretation, since we are not assuming that $A$ 
is a monoid.  We are saved by the fact that (since $k\geq 2$) either $b'_{ij}$
or $b''_{ij}$ is always the basepoint, and thus either $x(b'_{ij})$
or $y(b''_{ij})$ is always the basepoint $*$ of $A$.  It therefore suffices to 
define $ a\cdot *=*\cdot a=a$ for all $a\in A$ to make equation (\ref{e6})
meaningful in this context.
Finally, we choose $e\in\O(0)$ and $\mu\in \O(2)$ to be the basepoints.
}
\end{example}

As background for our final example we recall that, if $R$ is a commutative
ring and $m\geq 2$, Pirashvili \cite{Pira} has defined a higher Hochschild 
homology complex $Y_\b(R)$ by letting $Y_n(R)$ be the tensor product
\[
\bigotimes_{B_n} R
\]
(where $B_\b$ is the simplicial set from example \ref{ex3}), with face and
degeneracy maps induced by those of $B_\b$; thus $Y_\b(R)$ is the tensor
product $B_\b \otimes R$ in the category of simplicial commutative rings, as
defined by Quillen \cite{QuilHA}.

We also recall \cite[Proposition VII.1.6]{EKMM} that $\wedge$ is the 
coproduct in the category of commutative $S$-algebras. It follows that if $R$ 
is an $S$-algebra, $A$ and $B$ any sets and 
$f: A\rightarrow B$ any map then $f$ induces a map
\[
f_*:\bigwedge_A R \rightarrow \bigwedge_B R
\]
of $S$-algebras.

\begin{example}{\bf (The Pirashvili cohomology complex of a commutative ring 
spectrum.)}
{\rm 
Fix $k\geq 2$. Let $R$ be a commutative $S$-algebra, and define the $k$-th
Pirashvili homology complex $Y_\b(R)$ by letting
\[
Y_n(R) = \bigwedge_{B_n} R
\]
(where $B_\b$ is the simplicial set from example \ref{ex3})
with face and degeneracy maps induced by
those of $B_\b$.  The inclusion of the basepoint in $B_n$ induces a map
$R \rightarrow Y_n(R)$ of commutative $S$-algebras which makes $Y_n(R)$ an
$R$-module.  We define the $k$-th Pirashvili cohomology complex of $R$ to be
the cosimplicial $R$-module $X^\b(R)$ with
\[
X^n(R)=F_R(Y_n(R),R)
\]
(see \cite[Section III.6]{EKMM} for the definition of $F_R$), with coface and 
codegeneracy maps induced by the face and degeneracy maps of $Y_\b$.  Next we 
want to define an operad with multiplication $\O$ whose associated 
cosimplicial object is $X^\b$.  Let $\O(n)=X^n$. In order to define the 
operation $\circ_i$ we first observe that, since $k\geq 2$, the map
\[
B_{n+j-1}\rightarrow B_n \times B_j
\]
which takes $b$ to $(b'_{ij},b''_{ij})$ factors through the wedge to give a map
\[
B_{n+j-1}\rightarrow B_n\vee B_j
\]
and this in turn induces a map
\[
f_i:Y_{n+j-1}(R)=\bigwedge_{B_{n+j-1}} R \rightarrow \bigwedge_{B_n\vee B_j} R
\]
We next observe that there is a natural isomorphism
\[
\bigwedge_{B_n\vee B_j} R \cong (\bigwedge_{B_n} R)\wedge_R (\bigwedge_{B_j} R)
= Y_n(R)\wedge_R Y_j(R)
\]
We can therefore define $\circ_i$ to be the composite
\begin{eqnarray*}
X^n\wedge X^j &=&
F_R(Y_n(R),R)\wedge F_R(Y_j(R),R) \labarrow{\wedge_R}
F_R(Y_n(R)\wedge_R Y_j(R),R\wedge_R R) \\
&\cong&
F_R(\bigwedge_{B_n\vee B_j} R,R) \labarrow{f_i^*}
F_R(Y_{n+j-1}(R),R)=X^{n+j-1}
\end{eqnarray*}
It remains to specify the maps $e:S\rightarrow \O(0)$ and $\mu:S\rightarrow
\O(2)$.  We note that $\O(0)$ is just $R$, so we can let $e$ be the unit map of
$R$.  If $k>2$ then $\O(2)$ is also $R$, and we can let $\mu$ be the unit map.
If $k=2$ then $\O(2)$ is isomorphic to $F_S(R,R)$, and we let $\mu$ be the 
 adjoint of the identity map.
}
\end{example}

\newpage

\section{The spaces of the operad $\C'$.}
\label{sec3}

Recall that $\C_2$ denotes the operad of little 2-cubes.  The zeroth space of
$\C_2$ is a point.  We can define a suboperad $\C_2'$ of $\C_2$ by letting
$\C_2'(n)=\C_2(n)$ for all $n>0$ and letting $\C_2'(0)$ be the empty set.
All of the information about $\C_2$ and its operad structure is contained in
$\C_2'$ {\it except} for the so-called ``degeneracy'' maps
\[
\circ_k:\C_2(n)\times\C_2(0)\rightarrow\C_2(n-1)
\]
(the corresponding maps for $\C_2'$ have the empty set as their domain).

In this section and the next we define an operad $\C'$ which will turn
out to be equivalent to $\C_2'$.  In section \ref{sec5} we will define the
operad $\C$ of theorem \ref{main} by adding one more ingredient to the
definition of $\C'$.

\subsection{The spaces $\P(n)$ and $\F(n)$.} \label{subseca}

The $n$-th space $\C'(n)$ is empty for $n=0$ and for $n>0$ it is the Cartesian 
product of a space $\F(n)$ and a contractible space $\P(n)$.  (The idea here 
is that the spaces $\F(n)$ by themselves have a structure like that of an 
operad but with a composition operation $\gamma$ which is only 
$A_\infty$-associative.  By combining the $\F(n)$ with the $\P(n)$ we get a 
true operad).

The space $\P(n)$ is easy to define: it is the empty set for $n=0$ and for
$n>0$ it is the set of $n$-tuples $(p_1,\ldots,p_n)$ of positive real 
numbers that add up to 1.

The structure of the space $\F(n)$ is related to that of the algebraic operad
$\H$ defined in the introduction.
We will describe it in several stages: first we define the
set that indexes the cells of $\F(n)$, then we define the cells themselves, and
finally we define the attaching maps.

\subsection{The set of ``formulas of type $n$.''}
\label{subsecz}

The cells of $\F(n)$ will be indexed by what we call ``formulas of type $n$.'' 
Specifically, a formula of type $n$ is a symbol constructed
from the numbers $1,2,\ldots n$ (with no repetitions) by formal 
cup products and formal operad compositions. (The formulas of type $n$ 
correspond to certain multilinear operations in the Hochschild complex.  The
formal cup products correspond to cup products in the Hochschild complex, and
the formal operad compositions correspond to brace operations.)
For ease of
notation we write the formal operad composition $\gamma(x_0;x_1,\ldots,x_n)$ as
$x_0(x_1,\ldots,x_n)$.

Some examples of formulas of type 5 are: $2(1(3,5),4)$, 
$1(2\smallsmile 4, 3(5))$, $4(2,3)\smallsmile 1(5)$.

Here is a more precise definition.  Define a sequence of sets $\{S_k\}$ for
$k\geq 0$ as follows: 
$S_k$ is the set of positive integers when $k\neq 2$ and $S_2$ is the set
of positive integers with the symbol $\mu$ adjoined.  Take the free operad
generated by the sequence of sets $S_k$ and impose the relation (\ref{e3})
from Section \ref{sec2}.  Let the resulting operad be denoted by $\J'$ (the
prime refers to the fact that we are engaged in defining $\C'$).  
Then the set of formulas of type $n$ is precisely the set of elements in
$\J'(0)$ which contain each of the symbols $1,\ldots,n$ exactly once.

\begin{remark}
\label{r1}
{\rm
For clarity we should point out that, in this way of defining the formulas of
type $n$,
the cup products are denoted by $\mu$ instead of $\smallsmile$, so that
for example the formula $1(2\smallsmile 4, 3(5))$ would be written
$1(\mu(2,4),3(5))$ and the formula $4(2,3)\smallsmile 1(5)$ would be written
$\mu(4(2,3),1(5))$.
}
\end{remark}

For our later definitions it will be convenient to define the {\it valence} of
an integer in $S_k$ to be the number $k$.  Thus in a formula of type $n$ the
valence of a symbol $1,\ldots,n$ is the number of formal inputs for that
symbol.  For example, in the formula 
\[
3(2,4(5,6),1,7)\smallsmile 8(9)
\]
the valence of the symbol $3$ is 4, the valence of the symbol $4$ is 2, the 
valence of the symbol 8 is 1 and the valences of $1,2,5,6,7,9$ are each 0.

We will write $v(i)$ for the valence of $i$.

\subsection{The cells of $\F(n)$.} \label{subseca1}

To each formula $f$ of type $n$ we associate the product of simplices
\[
\prod_{i=1}^n \Delta^{v(i)}
\]
where $v(i)$ is the valence of the symbol $i$.
This product of simplices is a cell of $\F(n)$. We will frequently denote 
the cell associated to a formula $f$ by $\F_f$.

Here are some examples of cells of $\F(5)$:
\[
\F_{2(1(3,5),4)}= \Delta^2\times \Delta^2 \times \Delta^0 \times 
\Delta^0\times \Delta^0,
\]
\[
\F_{1(2\smallsmile 4, 3(5))} 
=\Delta^2\times \Delta^0 \times \Delta^1
\times \Delta^0\times \Delta^0,
\]
and 
\[
\F_{4(2,3)\smallsmile 1(5)}=
\Delta^1\times 
\Delta^0\times \Delta^0 \times
\Delta^2 \times \Delta^0.
\]

\subsection{Identifications Along The Boundary.} \label{subseca2}

The boundary of each cell $\F_f$ is a union of other cells which are selected 
by a rule modeled on the relation (\ref{tu9}) given in the introduction.  
Before spelling this out we give some examples of what the rule will say.

\bigskip

\noindent
{\bf Example.} 
The cell $\F_{1(2)}$ is an interval whose initial point is
identified with $\F_{1\smallsmile 2}$ and whose terminal point is identified
with $\F_{2\smallsmile 1}$. 
%(Compare this with the formulas 
%\[
%(d^0 r)(r_0)=r_0\cdot r, \qquad (d^1 r)(r_0)= r\cdot r_0
%\]
%which hold when $r$ is
%a 0-cochain in the Hochschild cohomology complex).
Similarly, the cell $\F_{2(1)}$ is an interval whose initial point is
identified with $\F_{2\smallsmile 1}$ and whose terminal point is identified
with $\F_{1\smallsmile 2}$. Since the only cells of $\F(2)$ are $\F_{1(2)}$,
$\F_{2(1)}$, $F_{1\smallsmile 2}$, and $\F_{2\smallsmile 1}$, we conclude that
the space $\F(2)$ is homeomorphic to a circle.

\bigskip

\noindent
{\bf Example.} $\F_{1(2,3)}$ is a triangle whose three 
faces
are identified with $\F_{2\smallsmile 1(3)}$, $\F_{1(2\smallsmile 3)}$ and
$\F_{1(2)\smallsmile 3}$.  
%(Compare this with the formulas
%\[
%(d^0 x)(r_0,r_1)=r_0 x(r_1), \qquad
%(d^1 x)(r_0,r_1)= x(r_0 r_1), \qquad
%(d^2 x)(r_0,r_1)= x(r_0) r_1
%\]
%which hold when $x$ is a 1-cochain in the
%Hochschild cohomology complex).

\bigskip

\noindent
{\bf Example.} $\F_{1(2(3))}$ is a square whose four faces are
identified with
$\F_{2(3)\smallsmile 1}$,
$\F_{1\smallsmile 2(3)}$,
$\F_{1(3\smallsmile 2)}$, and
$\F_{1(2\smallsmile 3)}$.

\bigskip

Next we give a formal description of the rule which underlies these 
examples.  We will use the operad $\J'$ described in
subsection \ref{subseca}.

Let $f$ be a formula of type $n$ and let  $i$ be an integer with $1\leq i \leq
n$.  Let $k$ be the valence of $i$ in $f$ and let $j$ be an 
integer with $0\leq j\leq k$. If $k\geq 1$ we define 
$\partial_{ij} f$ to be the formula obtained from $f$ by replacing the symbol 
$i$ by 
\[
\left\{
\begin{array}{ll}
\mu\circ_2 i_{k-1}
& \mbox{if $j=0$} \\
i_{k-1}\circ_j \mu & \mbox{if $0< j < k$} \\
\mu\circ_1 i_{k-1} & \mbox{if $j=k$}
\end{array}
\right.
\]
Here we are writing $i_k$ for the copy of $i$ in $S_k$; also see Remark 
\ref{r1}. 
Of course, this definition is motivated by equation (\ref{e1}).
Here are some examples: 
\[
\partial_{10}\, 1(2)=2\smallsmile 1 \ \mbox{and} \ 
\partial_{11}\, 1(2)=1\smallsmile 2 
\]
\[
\partial_{10}\, 1(2,3)=2\smallsmile 1(3), \quad \partial_{11}\,
1(2,3)=1(2\smallsmile 3), \ \mbox{and}\ \partial_{12}\,1(2,3)=1(2)\smallsmile 3
\]
\[
\partial_{10}\, 1(2(3))=2(3)\smallsmile 1, \quad \partial_{11}\,
1(2(3))=1\smallsmile 2(3)
\]
\[
\partial_{20}\, 1(2(3))=1(3\smallsmile 2), \quad
\partial_{21}\, 1(2(3))=1(2\smallsmile 3)
\]

Now let ${\cal I}_n$ be the set of formulas of type $n$.  We can define a 
partial 
ordering on ${\cal I}_n$ as follows: given a formula $f$ of type $n$, the 
maximal 
elements in the set 
$\{\, g\,|\, g < f \,\}$ are $\{\, \partial_{ij} f \,|\, 1\leq i\leq n, 0\leq
j\leq v(i)\,\}$.

Let us now think of the partially ordered set
${\cal I}_n$
as a category in the usual way.  We can make the function $\F$ defined in the 
previous subsection 
into a functor from ${\cal I}_n$ to the category of topological 
spaces by taking the arrow $\partial_{ij} f\rightarrow f$ to the map
\[
1\times\cdots\times d_j\times\cdots\times 1:
\prod_{k=1}^{i-1} \Delta^{v(k)}\times \Delta^{v(i)-1}\times\prod_{k=i+1}^n
\Delta^{v(k)}
\rightarrow 
\prod_{k=1}^n \Delta^{v(k)} 
%\F_{\partial_{ij}f}\rightarrow\F_f,
\]
where $d_j$ is the $j$-th face map 
\[
\Delta^{v(i)-1}\rightarrow\Delta^{v(i)}
\]

Now we can give a formal definition of the space $\F(n)$:

\begin{definition}
$\F(n)=\mbox{colim}_{{\cal I}_n} \ \F$
\end{definition}

This definition implies that the boundary of a cell $\F_f$ is the union of the
cells $\F_{\partial_{ij} f}$.

\newpage

\section{The operad structure of $\C'$ and the action of $\C'$ on ${\mbox{\rm Tot}}(X^\b)$.}
\label{sec4}

Our next goal is to do two things: to define the operad structure of $\C'$ and
to define the action maps
\[
\theta_n:\C'(n)\times ({\mbox{\rm Tot}}(X^\b))^n \rightarrow {\mbox{\rm Tot}}(X^\b)
\]
when $X^\b$ is a cosimplicial space arising from an operad with multiplication
(the case where $X^\b$ is a cosimplicial $S$-module requires only routine
changes and we will not discuss it separately).  
It turns out that the second of these 
tasks is easier and provides helpful background for the first, so we begin 
with it and defer the definition of the operad structure to the end of this
section.

We will construct $\theta_n$ as a map of sets and leave it to the reader to
verify continuity.

Recall that $\C'(n)=\P(n)\times\F(n)$.  In order to construct $\theta_n$ it
suffices to construct, for each choice of $\p=(p_1,\ldots,p_n)\in P(n)$ and of
a cell $\F_f$ of $\F(n)$,
a map 
\[
\theta_{\p,f}:\F_f\times ({\mbox{\rm Tot}}(X^\b))^n \rightarrow {\mbox{\rm Tot}}(X^\b)
\]
Next recall that 
\[
{\mbox{\rm Tot}}(X^\b)=\Hom(\Delta^\b,X^\b),
\]
where $\Hom$ denotes morphisms of cosimplicial spaces, so in order to construct
$\theta_n$ it suffices to construct, for each $k\geq 0$, a suitable map
\[
\tilde{\theta}_{\p,f,k}:
\Delta^k\times\F_f
\times 
({\mbox{\rm Tot}}(X^\b))^n
\rightarrow 
X^k.
\]
We do the case $k=0$ first, since this illustrates the general idea in a simple
situation.  

\begin{example}[The case $k=0$.] 
\label{exam1}
{\rm
Given a formula $f$ of type $n$ and elements $x_i\in X^{v(i)}$ for
$1\leq i\leq n$, we can define an element $\bar{f}(x_1,\ldots,x_n)\in X^0$ by
replacing each symbol $i$ in $f$ by $x_i$ and then interpreting the formal
compositions and cup products in $f$ to be genuine compositions and cup
products in $X^\b$.  For example, if $f=3(1(2\smallsmile 4),6(5))$ then
\[
\bar{f}(x_1,\ldots,x_6)=x_3(x_1(\mu(x_2,x_4)),x_6(x_5)).
\]
This process gives a map
\[
\bar{f}:\prod_{i=1}^n X^{v(i)} \rightarrow X^0.
\]
Now let $\s\in\F_f$, and ${\mathbf a}_i\in{\mbox{\rm Tot}}(X^\b)$ for 
$1\leq i\leq n$; thus $\s$ is an $n$-tuple $\s_1,\ldots,\s_n$ with 
$\s_i\in\Delta^{v(i)}$ and
each ${\mathbf a}_i$ is a sequence 
$\{\,a_{im}:\Delta^m\rightarrow 
X^m\,\}$.  We define
\[
\tilde{\theta}_{\p,f,0}(\s,{\mathbf a}_1,\ldots,{\mathbf a}_n)=
\bar{f}(a_{1v(1)}(\s_1),\ldots,
a_{nv(n)}(\s_n)).
\]
(Note that $\p$ plays no role in this definition, but it will have a role for
$k>0$).
}
\end{example}

In order to extend this idea to $k>0$ we need a way of decomposing the product
$\Delta^k\times \F_f$ as a union of products of 
the form $\prod_{i=1}^n \Delta^{k_i}$.  Here is an example which should 
help the reader to follow the general description. It may also be helpful to
compare this example to the relation (\ref{tu777}) for iterated brace products
given in the introduction.

\begin{example} 
\label{exam2}
{\rm
The following
picture shows how to define the map
\[
\tilde{\theta}_{\p,f,k}:
\Delta^k \times  \F_f 
\times ({\mbox{\rm Tot}}(X^\b))^n
\rightarrow
X^k
\]
for the case in which $n$ is 2, $f$ is $1(2)$, $k$ is 1, and $\p=(1/2,1/2)$ 
%is an arbitrary point of $\P(2)$
.
\medskip

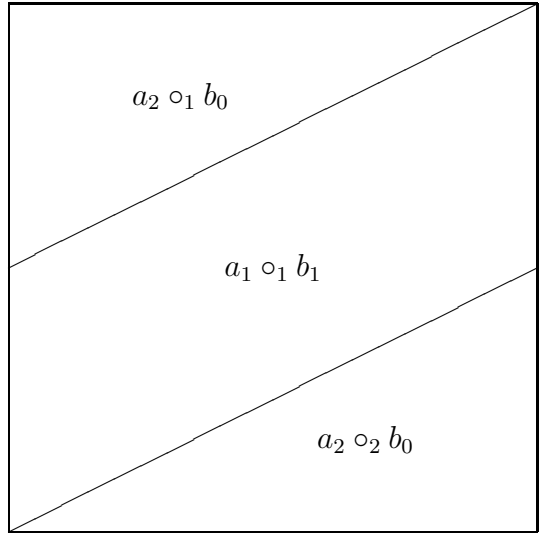
\hfill
\begin{picture}(200,200)

\put(0,0){\line(1,0){200}}
\put(0,0){\line(0,1){200}}
\put(0,200){\line(1,0){200}}
\put(200,0){\line(0,1){200}}

\put(0,0){\line(2,1){200}}
\put(0,100){\line(2,1){200}}

\put(0,95){\makebox(200,10){$a_1\circ_1 b_1$}}

\put(0,160){\makebox(130,10){$a_2\circ_1 b_0$}}
\put(70,30){\makebox(130,10){$a_2\circ_2 b_0$}}

%\put(-30,45){\makebox(30,10){$p_2$}}
%\put(200,45){\makebox(30,10){$p_1$}}

%\put(-30,145){\makebox(30,10){$p_1$}}
%\put(200,145){\makebox(30,10){$p_2$}}

\end{picture}
\hfill{\mbox{}}

\medskip
\noindent
Here ${\mathbf a}=\{a_m\}$ and ${\mathbf b}=\{b_m\}$ are points of
${\mbox{\rm Tot}}(X^\b)$, and the picture defines a map from
$\Delta^1\times \F_{1(2)}$ to $X^1$ (note: the $\F_{1(2)}$ coordinate
is the horizontal direction).  This picture is part of a homotopy from
${\mathbf b}\smallsmile {\mathbf a}$ to ${\mathbf a}\smallsmile
{\mathbf b}$, since the left- and right-hand edges are the projections
of ${\mathbf b}\smallsmile {\mathbf a}$ and ${\mathbf a}\smallsmile
{\mathbf b}$ on $\Hom(\Delta^1,X^1)$.  }
\end{example}

Notice that each vertical cross-section of the picture in example \ref{exam2}
is a 3-fold prismatic subdivision of $\Delta^1$ (except at the endpoints where
it degenerates to a 2-fold prismatic subdivision).  We therefore think of the 
picture in example \ref{exam2} as a ``fiberwise prismatic subdivision'' of 
$\Delta^1 \times \F_{1(2)}$.  Our next goal is to define the 
fiberwise prismatic subdivision of $\Delta^k\times \F_f$ in general, and for 
this we first need to define the ``thickenings'' of a formula $f$.

\subsection{The thickenings of a formula.}

\label{subsecb}

The cells of the fiberwise prismatic subdivision of $\Delta^k\times F_f$ will 
correspond to a collection of formulas called the $k$-thickenings of $f$.
To thicken a formula means to add extra inputs (denoted by the symbol $\id$) to 
the formula (we will give a formal definition in a moment).  The new formula 
is called a $k$-thickening if there are $k$ copies of the symbol $\id$.

Before giving the formal definition we give an example.
The 1-thickenings of 1(2) are $1(\id,2)$, $1(2,\id)$, 
$1(2(\id))$, $1(\id \smallsmile 2)$, $1(2\smallsmile \id)$, $\id\smallsmile
1(2))$ and $1(2)\smallsmile \id$,
and in Example \ref{exam2} these correspond to the cells which are labeled by
$a_2\circ_2 b_0$, $a_2\circ_1 b_0$, $a_1\circ_1 b_1$, and to the two diagonal
and the two horizontal edges, respectively.

For the formal definition, we use the operad $\J'$ defined in subsection
\ref{subseca}.  Let us define a formula of type $(n,k)$ to be an element of 
$\J'(k)$ which contains each of the symbols $1,\ldots,n$ exactly once (note 
that such a formula is forced to contain the symbol $\id$ exactly $k$ times).  
There is a reduction map
\[
\rho:\J'(k)\rightarrow \J'(0)
\]
which removes each copy of $\id$, together with the typographical symbols 
that are immediately adjacent to it (these may be commas, parentheses, or
$\smallsmile$'s).  For example,
$\rho$ takes $1(2\smallsmile \id,\id,3(4(\id)))$ to $1(2,3(4))$.

\begin{definition}
{\rm
If $f$ is a formula of type $n$, the $k$-thickenings of $f$ are the formulas
$g$ of type $(n,k)$ for which $\rho(g)=f$.
}
\end{definition}

\subsection{Fiberwise prismatic subdivision.}
\label{subsecbb} 

Given a formula $f$ of type $n$, let $\I_{f,k}$ be the set of
$k$-thickenings of $f$.  We can define a function $\F$ from $\I_{f,k}$
to the set of topological spaces exactly as in subsection
\ref{subseca1}, that is,
\[
\F_g=\prod_{i=1}^n \Delta^{v(i)}
\]
where $g$ is in $\I_{f,k}$ and $v(i)$ denotes the valence of $i$ in $g$.

For our present purposes we need to generalize this, using simplices
of variable sizes as in Section \ref{sec1} : given
$\q=(q_1,\ldots,q_n)\in ({\mathbb R}_{\geq 0})^n$ and $g\in \I_{f,k}$ we
define
\[
\F_{g,\q}=\prod_{i=1}^n \Delta^{v(i)}_{q_i};
\]
this gives a function 
$\F_{-,\q}$ from $\I_{f,k}$ to the set of topological spaces for each choice of
$\q$.

Next observe that
we can make $\I_{f,k}$ into a partially ordered set, and $\F_{-,\q}$ into a
functor from $\I_{f,k}$ to the category of topological spaces, exactly as we
did with $\I_n$ in subsection \ref{subseca2}. 

\begin{definition}
{\rm
$\F_{f,k,\q}$ is the space ${\mbox{colim}}_{{\cal I}_{f,k}} \ \F_{-,\q}$.
}
\end{definition}

(In example \ref{exam2}, $\F_{1(2),1,(1,1)}$ is the union of two triangles 
and a square, with the edge identifications indicated in the picture there.)

The fiberwise prismatic subdivision is a certain map
\[
\sigma:\F_{f,k,\q}\rightarrow\Delta^k_{q_1+\cdots+q_n}\times \F_{f,\q}
\]
which will turn out to be a homeomorphism.
In order to describe it we need some terminology.

Let $g$ be an element 
of $\I_{f,k}$.  Then $g$ is a character-string consisting of the integers
1,\dots,$n$ together with $k$ copies of the symbol $\id$ and some commas and 
parentheses.  If $i\in\{\,1,\ldots,n\,\}$ and $v(i)>0$ then there are 
$v(i)-1$ commas, a left parenthesis, and a right parenthesis which {\it 
belong} to $i$ (to be precise, the characters that belong to $i$ are the left 
and right parentheses that enclose the inputs of $i$ and the commas that 
separate them).  If $i$ has valence 0 then we say that $i$ belongs to itself. 
Note that in either case there are $v(i)+1 $ characters which 
belong to $i$. We refer to the characters of $g$ which belong to some $i$ as 
the {\it eligible} characters of $g$; thus the eligible characters are all of 
the commas, all of the parentheses, and the $i$ with valence 0.  

Next let $\s$ be an element of $\F_{g,\q}$.  Now $\F_{g,\q}$ is a product of
simplices indexed by 1,\dots,n, and we write $\s_i$ for the projection of 
$\s$ on the $i$-th factor; thus $\s_i\in\Delta^{v(i)}_{q_i}$.  Let us fix 
$i$ temporarily. As an element of $\Delta^{v(i)}_{q_i}$, $\s_i$ has $v(i)+1$ 
coordinates, and we denote these by $s_{ij}$.  Since there are $v(i)+1$ 
characters belonging to $i$ in $g$, we can match the numbers $s_{ij}$ with 
the characters that belong to $i$ (going from left to right as $j$ 
increases), and letting $i$ vary we get a one-to-one correspondence between 
the eligible characters of $g$ and the numbers $s_{ij}$, $1\leq i\leq n$, 
$0\leq j\leq v(i)$.  Here is an example, for the formula
$g=3(2(1,\id),4,\id,\id,5(6))$: 
\[
\begin{array}{cccccccccccccccccccc}
3 & (     & 2 & (    & 1    & ,    & \id &    ) &    , &    4 & , &\id&, &\id&,    
& 5 &    
(     &    6 &    ) &    ) \\
  & s_{30}&   &s_{20}&s_{10}&s_{21}&     &s_{22}&s_{31}&s_{40}&s_{32} &&
s_{33}&&
s_{34}
&   &
s_{50}&s_{60}&s_{51} &s_{35}
\end{array}
\]
We will refer to this diagram as the {\it tableau} of $g$. The tableau of $g$ 
has two lines, with the first being the string of characters of $g$ and the 
second the symbols $s_{ij}$ in the order we have prescribed.

Now we are ready to describe the fiberwise prismatic subdivision map 
\[
\sigma:\F_{f,k,\q}\rightarrow\Delta^k_{q_1+\cdots+q_n}\times \F_{f,\q}
\]
We write $\sigma_2$ for the projection of $\sigma$ on the second factor: 
\[
\sigma_2:\F_{f,k,\q}\rightarrow\F_{f,\q}
\]
Since $\F_{f,\q}$ is itself a product, it suffices to 
define the projection of $\sigma_2$ on the $i$-th factor of $\F_{f,\q}$; we 
denote this projection by $\sigma_{2i}$.  
Let $g\in \I_{f,k}$ and let
$\s=(\s_1,\ldots,\s_n)\in \F_{g,\q}$. 
In order to apply 
$\sigma_{2i}$ to $\s$, we select the coordinates $s_{ij}$ from the second line
of the tableau of $g$ and then insert $+$ signs whenever the symbol $\id$
occurs in the first line. 
(The motivation for this is that $f$ is
obtained from $g$ by collapsing out the symbols $\id$ that occur in $g$.)
For example, in the tableau given above we have
\[
\begin{array}{ccl}
\sigma_{21}(\s)&=&s_{10} \\
\sigma_{22}(\s)&=& (s_{20},s_{21}+s_{22}) \\
\sigma_{23}(\s)&=& (s_{30},s_{31},s_{32}+s_{33}+s_{34},s_{35}) \\
\sigma_{24}(\s)&=& s_{40} \\
\sigma_{25}(\s)&=& (s_{50},s_{51}) \\
\sigma_{26}(\s)&=& s_{60}
\end{array}
\]
(Note that the restriction of $\sigma_2$ to $\F_{g,\q}$ is a Cartesian 
product of iterated degeneracies, with the $\id$ symbols telling what
degeneracies are to be used.)

Next we define the projection of $\sigma$ on the first factor, which we denote
by $\sigma_1$:
\[
\sigma_1:\F_{f,k,\q}\rightarrow\Delta^k_{q_1+\cdots+q_n}
\]
Again let
$g\in \I_{f,k}$ and let
$\s=(\s_1,\ldots,\s_n)\in \F_{g,\q}$. 
In order to apply $\sigma_1$ to $\s$, we insert $+$ signs in the second line of
the tableau wherever the symbol $\id$ does {\it not} occur in the first 
line.  For example, in the tableau given above we have $k=3$ and 
\[
\sigma_1(\s)= (s_{30}+s_{20}+s_{10}+s_{21},s_{22}+s_{31}+s_{40}+s_{32},s_{33},
s_{34}+s_{50}+s_{60}+s_{51}+s_{35}).
\]

When these definitions are applied to example \ref{exam2}, the formulas that
result are:

\begin{itemize}
\item If $\s\in\F_{1(\id,2),\q}$ then 
$\sigma_1(\s)=(s_{10},s_{11}+s_{20}+s_{12})$ and 
$\sigma_2(\s)=((s_{10}+s_{11},s_{12}),s_{20})$.
\item If $\s\in\F_{1(2,\id),\q}$ then 
$\sigma_1(\s)=(s_{10}+s_{20}+s_{11},s_{12})$ and 
$\sigma_2(\s)=((s_{10},s_{11}+s_{12}),s_{20})$.
\item If $\s\in\F_{1(2(\id)),\q}$ then 
$\sigma_1(\s)=(s_{10}+s_{20},s_{21}+s_{11})$ and
$\sigma_2(\s)=((s_{10},s_{11}),s_{20}+s_{21})$.
\end{itemize}

Returning to the general situation, it is not difficult to see that the 
restriction of $\sigma_1$ to each fiber of $\sigma_2$ is a prismatic 
subdivision of $\Delta^k_{q_1+\cdots+q_n}$; this accounts for the name 
``fiberwise prismatic subdivision'' and also implies

\begin{proposition}
$\sigma$ is a homeomorphism.
\end{proposition}

Finally, it is convenient to introduce a variant of $\sigma$.  Let
$\p\in\P(n)$ and let $f$ be a formula of type $n$. For each $g\in\I_{f,k}$ 
define 
\[
\cdot\p:\F_g\rightarrow \F_{g,\p}
\]
to be the map which takes 
$(\s_1,\ldots,\s_n)\in\F_g$ to 
$(p_1\s_1,\ldots,p_n\s_n)\in\F_{g,\p}$.
These maps fit together to give a map
\[
\cdot\p:\F_{f,k}\rightarrow\F_{f,k,\p}
\]
(where, as the reader may have guessed, $\F_{f,k}$ is an abbreviation for
$\F_{f,k,(1,\ldots,1)}$).
Now define
\[
\sigma(\p):\F_{f,k}\rightarrow \Delta^k \times \F_f
\]
to be the composite
\[
\F_{f,k}\labarrow{\cdot\p}
\F_{f,k,\p} \labarrow{\sigma}
\Delta^k\times\F_{f,\p}
\labarrow{1\times (\cdot\p)^{-1}}
\Delta^k\times\F_f
\]
Note that the projection of $\sigma(\p)$ on the second factor $\F_f$ is the
same as the projection of $\sigma$ on $\F_f$ (that is, it is the 
map $\sigma_2$ defined above).

\subsection{Definition of the map $\tilde{\theta}_{\p,f,k}$}

We are now ready to define the map
\[
\tilde{\theta}_{\p,f,k}:
\Delta^k \times  \F_f 
\times ({\mbox{\rm Tot}}(X^\b))^n
\rightarrow
X^k;
\]
it is the composite
\[
\begin{array}{l}
\Delta^k \times  \F_f 
\times ({\mbox{\rm Tot}}(X^\b))^n
\labarrow{\sigma(\p)^{-1}\times 1}
\F_{f,k} \times ({\mbox{\rm Tot}}(X^\b))^n
\labarrow{\theta'_{f,k}}
X^k,
\end{array}
\]
where $\sigma(\p)$ was defined at the end of subsection \ref{subsecbb} and 
$\theta'_{f,k}$ is a map which we define next.

Recall that $\F_{f,k}$ is 
${\mbox{colim}}_{g\in{\cal I}_{f,k}} \ \F_{g}$.
It therefore suffices to define the restriction of $\theta'_{f,k}$ to $\F_g$;
we denote this restriction by
\[
\theta'_{g}:\F_{g}
\times ({\mbox{\rm Tot}}(X^\b))^n
\rightarrow X^k.
\]

So fix $g\in\I_{f,k}$. As in Example \ref{exam1}, if we are given 
elements $x_i\in X^{v(i)}$ for
$1\leq i\leq n$ (where $v(i)$ denotes the valence of $i$ in $g$) we can 
define an element $\bar{g}(x_1,\ldots,x_n)\in X^k$ by
replacing each symbol $i$ in $g$ by $x_i$ and then interpreting the formal
compositions and cup products in $g$ to be genuine compositions and cup
products in $X^\b$;  
this process gives a map
\[
\bar{g}:\prod_{i=1}^n X^{v(i)} \rightarrow X^k.
\]

Now given $\s\in\F_{g}$ and 
${\mathbf a}_i\in{\mbox{\rm Tot}}(X^\b)$ for $1\leq i\leq n$
we define
\[
\theta'_{g}(\s,{\mathbf a}_1,\ldots,{\mathbf a}_n)=
\bar{g}\left(a_{1v(1)}(\s_1),\ldots,
a_{nv(n)}(\s_n)\right),
\]
where $v(i)$ denotes the valence of $i$ in $g$.

This completes the definition of the maps $\tilde{\theta}_{\p,f,k}$.
It is straightforward to check that these maps fit together to give the map
$\theta_n:\C'(n)\times ({\mbox{\rm Tot}}(X^\b))^n \rightarrow {\mbox{\rm Tot}}(X^\b)$
that we set out to define.

\subsection{The operad structure of $\C'$.}
\label{subsecc}

The operad structure of $\C'$ is determined by the fact that we want $\theta_n$ 
to be an action of $\C'$.

First we describe the action of the symmetric group $\Sigma_n$ on 
$\C'(n)$ (we depart slightly from \cite{MayG} by having the symmetric group 
act on the left instead of the right.) Let $\tau\in\Sigma_n$.  
The action of $\tau$ on the $\P(n)$ factor permutes the
coordinates $(p_1,\ldots,p_n)$.
Given a formula $f$ of type $n$, we define $\tau(f)$ to
be the formula obtained from $f$ by replacing each symbol $i$ in $f$ by
$\tau(i)$.  There is an evident homeomorphism 
$\tau:\F_f\rightarrow\F_{\tau(f)}$
which permutes the coordinates, and passing to the colimit over $f$ 
we get a homeomorphism $\tau:\F(n)\rightarrow\F(n)$.

Next we will describe the operad composition $\circ_k$.  We begin by giving the 
collection of spaces $\P(n)$ an operad composition: we define
\[
\circ_k:\P(n)\times \P(j)\rightarrow \P(n+j-1)
\] 
to be the map that takes the pair 
\[
(p_1,\ldots,p_n), (q_1,\ldots,q_j)
\]
to
\[
(p_1,\ldots,p_{k-1}, p_kq_1,\ldots,p_k q_j,p_{k+1},\ldots,p_n).
\]
(If we think of an element of $\P(n)$ as a collection of little intervals whose
lengths add up to 1 then $\P(n)$ is a suboperad of the little intervals 
operad $\C'_1$.)

Now we want to define
\[
\circ_k:
(\P(n)\times\F(n))
\times
(\P(j)\times\F(j))
\rightarrow
(\P(n+j-1)\times\F(n+j-1))
\]
for $1\leq i \leq n$.
The projection of $\circ_k$ on the $\P(n+j-1)$ factor is defined to be the
composite
\[
(\P(n)\times\F(n))
\times
(\P(j)\times\F(j))
\rightarrow
\P(n)\times\P(j)
\rightarrow
(\P(n+j-1),
\]
where the first map is the projection and the second is the $\circ_k$ operation
for the operad $\P$.
The projection of $\circ_k$ on the $\F(n+j-1)$ factor is defined to be the 
composite
\[
(\P(n)\times\F(n))
\times
(\P(j)\times\F(j))
\rightarrow
\P(j)\times\F(n)\times\F(j) 
\labarrow{c_{i,n,j}}
\F(n+j-1),
\]
where the first map is the projection and the map $c_{k,n,j}$ will be defined 
next.
(This means that the operad structure of $\C'$ is like a semidirect 
product of $\P$ with $\F$, except that $\F$ is not an operad.)

In order to construct $c_{k,n,j}$ it
suffices to construct, for each choice of $\p=(p_1,\ldots,p_n)\in P(n)$, 
each formula $f$ of type $n$, and each formula $f'$ of type $j$, a map 
\[
c_{k,\p,f,f'}:\F_f\times\F_{f'} \rightarrow \F(n+j-1)
\]
Let us write $\F_f^{\neq k}$ for 
\[
\prod_{i\neq k} \Delta^{v(i)};
\]
thus $\F_f=\F_f^{\neq k}\times \Delta^{v(k)}$, where $v(k)$ denotes the valence
of $k$ in $f$.
We define $c_{k,\p,f,f'}$ to be the composite
\[
\F_f\times\F_{f'}=\F_f^{\neq k}\times \Delta^{v(k)} \times \F_{f'}
\labarrow{1\times \sigma(\p)^{-1}}
\F_f^{\neq k}\times\F_{f',v(k)}
\labarrow{c'_{k,f,f'}}
\F_{n+j-1},
\]
where $c'_{k,f,f'}$ remains to be defined.

Of course, it suffices to define the restriction of $c'_{k,f,f'}$ to 
$\F_f^{\neq k}\times \F_{g'}$ where $g'$ is a $v(k)$-thickening of $f'$;
we denote this restriction by
\[
c'_{k,f,g'}:\F_f^{\neq k}\times\F_{g'} \rightarrow \F_{n+j-1}
\]

Next we define a formula $f*_k g'$ of type $n+j-1$ by ``substituting $g'$ for
$k$,'' as follows: we replace the symbols $k+1,\ldots,n$ in $f$ by
$k+j,\ldots,n+j-1$ respectively, replace the symbols $1,\ldots,j$ in $g'$ by 
$k,\ldots,k+j-1$ respectively, replace the symbols $\id$ in $g'$ by the entries
of $k$ in $f$, and then replace $k$ by $g'$.
For example, if $k=1$, $f=1(2(3),4\smallsmile 5, 6(7,8)$ and 
$g'=1(2,\id,3,\id,\id)$ then $f*_k g'=1(2,4(5),3,6\smallsmile 7, 8(9,10))$.

Observe that 
$\F_f^{\neq k}\times \F_{g'}
=\F_{f*_k g'}$. Since $f*_k g'$ is a formula of type $n+j-1$,
we can define $c'_{k,f,g'}$ to be the composite
\[
\F_f^{\neq k}\times \F_{g'}
=\F_{f*_k g'} \subset \F_{n+j-1},
\]
and this completes the definition of the operations $\circ_k$ 
for $\C'$.

The fact that the operations $\circ_k$ so defined make $\C'$ into an operad, and
the fact that the maps $\theta_n$ defined earlier give an action of this operad
on ${\mbox{\rm Tot}}(X^\b)$ when $X^\b$ is the cosimplicial space associated to an operad
with multiplication, are both easy consequences of an associativity property of
the fiberwise prismatic subdivision which will be stated and proved in the next
subsection.  First we pause to give two examples to illustrate the definition
of $\circ_k$.

\begin{example}
{\rm
Let $f=f'=1(2)$, let $\p=(p_1,p_2)$, and let $k=2$.
Then $v(k)=0$, so there is only one $g'$ (which is equal to $f'$) and
$c_{k,\p,f,f'}$ is the composite
\[
\F_f\times\F_{f'} =\F_{1(2(3))}\subset\F(3)
\]
}
\end{example}

\begin{example}
{\rm
Again let $f=f'=1(2)$, $\p=(p_1,p_2)$, but now let $k=1$.
Then $v(k)=1$, and the top-dimensional 1-thickenings $g'$ of $f'$ are
$1(\id,2)$, $1(2,\id)$, and $1(2(\id)$.
Now $c_{k,\p,f,f'}$ is the map of $\F_f\times\F_{f'}$ into $\F(3)$ indicated in
the following picture.

\bigskip
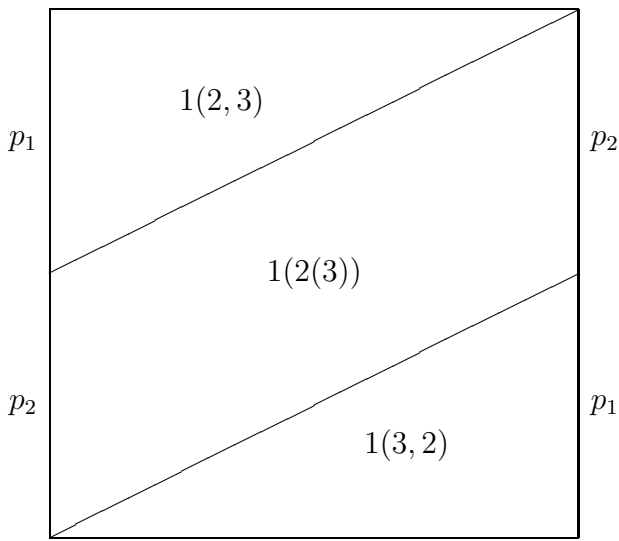
\begin{picture}(200,200)

\put(0,0){\line(1,0){200}}
\put(0,0){\line(0,1){200}}
\put(0,200){\line(1,0){200}}
\put(200,0){\line(0,1){200}}

\put(0,0){\line(2,1){200}}
\put(0,100){\line(2,1){200}}

\put(0,95){\makebox(200,10){$1(2(3))$}}

\put(0,160){\makebox(130,10){$1(2,3)$}}
\put(70,30){\makebox(130,10){$1(3,2)$}}

\put(-20,45){\makebox(20,10){$p_2$}}
\put(200,45){\makebox(20,10){$p_1$}}

\put(-20,145){\makebox(20,10){$p_1$}}
\put(200,145){\makebox(20,10){$p_2$}}

\end{picture}

}
\end{example}

\subsection{An associativity property of the fiberwise prismatic subdivision.}
\label{subsecd}

Given $k$, $n$, $j$ and $l$,
let $f$ be a formula of type $n$, $f'$ a formula of type $j$, $g$ 
an $l$-thickening of $f$, 
$g'$ a $v_f(k)$ thickening of $f'$
(where, of course,
$v_f(k)$ denotes the valence of $k$ in $f$), and
$h'$ a $v_g(k)$-thickening of $f'$. 
Note that $g*_k h'$ is an $l$-thickening of $f*_k g'$.
Let $\p\in\P(n)$ and let $\p'\in\P(j)$.

It is straightforward to check from the definitions that the following diagram
commutes.  

\[
\begin{array}{ccc}
\F_{g*_k h'} &=& 
\F_{g}^{\neq k} \times \F_{h'} \\
&& \\
\llap{$\sigma(\p\circ_k \p')$}\downarrow   && \downarrow \rlap{$1\times 
\sigma(\p')$} \\
&& \\
\Delta^l \times \F_{f*_k g'} && 
\F_{g}^{\neq k} \times \Delta^{v_g(k)} \times \F_{f'} \\
&& \\
\llap{$=$} \downarrow && \downarrow \rlap{$=$} \\
&& \\
\Delta^l \times \F_{f}^{\neq k}\times \F_{g'} && 
\F_{g} \times \F_{f'} \\
&& \\
\llap{$1\times1\times \sigma(\p')$} \downarrow & & \downarrow \rlap{$\sigma(\p) 
\times 1$}\\
&& \\
\Delta^l \times \F_{f}^{\neq k}\times \Delta^{v_f(k)} \times
\F_{f'} & = & 
\Delta^l\times \F_{f} \times \F_{f'} 
\end{array}
\]

\newpage

\section{The operad $\C$ and its action on ${\mbox{\rm Tot}}(X^\b)$.}
\label{sec5}

In this section we modify the definition of $\C'$ to obtain the operad $\C$ of
theorem \ref{main}.  

As we have already mentioned, the reason we need to modify $\C'$ is that it
doesn't model the degeneracy maps 
\[
\circ_k:\C_2(n)\times\C_2(0)\rightarrow \C_2(n-1)
\]
of the little 2-cubes operad.
Here is another way to describe the difficulty.  If $X^\b$ is as in Theorem
\ref{main}, the operad $\C'$ maps to a suboperad of the endomorphism operad of
Tot($X^\b$), but this suboperad is not closed under the degeneracy maps of
this endomorphism operad (these are the maps that insert one or more copies of
the basepoint).  We will remedy this defect by adjoining new points to 
each $\C'(n-1)$ which will serve as the degeneracies of
points in $\C'(n)$.

Before proceeding we need to specify the basepoint 
of ${\mbox{\rm Tot}}(X^\b)$.
The definition of operad with multiplication provides a special
element $e\in X^0$.  The conditions (\ref{e3})
and (\ref{e4}) of Section \ref{sec2} imply that, for each $n$, all $n$-fold 
iterated cofaces of $e$ are equal, so there is a unique cosimplicial map from 
$\Delta^\b$ to $X^\b$ which is constant on each $\Delta^n$ and takes 
$\Delta^0$ to $e$.   This cosimplicial map is the basepoint of ${\mbox{\rm Tot}}(X^\b)$; we 
will denote it by $\bar{e}$ from now on.

\begin{remark}\label{hi}
{\rm 
For later use we describe $\bar{e}$ more explicitly.  The projection of
$\bar{e}$ on $\Hom(\Delta^m,X^m)$ will be denoted by $\bar{e}_m$.  It is a 
constant map whose image is $e\in X^0$
for $m=0$, $\id\in X^1$ for $m=1$, $\mu\in X^2$ for $m=2$, and for $m\geq 3$ 
its image is the ``iterated multiplication'' $\mu\circ_1 (\mu\circ_1 ( \cdots 
\mu)\cdots)$ in $X^m$.
}
\end{remark}

In order to define $\C$
we follow the general outline of sections \ref{sec3} and \ref{sec4}.  First we
need to define an ``indexing'' operad $\J$ analogous to the operad $\J'$
defined in section \ref{sec3}. Define a sequence of sets $\{S_m\}$ for
$m\geq 0$ as follows:
$S_m$ is the set of positive integers with the symbol $\varepsilon$ adjoined
when $m\neq 2$ and $S_2$ is the set
of positive integers with the symbols $\varepsilon$ and $\mu$ adjoined.  Take 
the free operad
generated by the sequence of sets $S_m$ and impose the relation (\ref{e3})
from Section \ref{sec2}.  Let the resulting operad be denoted by $\J$. 

Let us write $\I_n^m$ for the set of elements of $\J(0)$ which contain 
each of the symbols $1,\ldots,n$ exactly once and the symbol $\varepsilon$ 
exactly $m$ times.  

If $f\in\I_n^k$, we define the valence of each of the symbols $1,\ldots,n$ and
of each of the copies of $\varepsilon$ in the usual way, as the number of 
entries of the symbol in question.  We define $\F_f$ to be
\[
\prod_{i=1}^n \Delta^{v(i)} \times 
\prod_{\varepsilon\in f}\Delta^{v(\varepsilon)}
\]
where the second product is indexed by the copies of $\varepsilon$ in $f$.

We can define a partial ordering on $\I_n^m$ and make $\F$ into a functor
from $\I_n^m$ to the category of topological spaces exactly as in subsection 
\ref{subseca2}.  We will denote the space 
$\mbox{colim}_{{\cal I}_n^m} \ \F$ by $\F(n,m)$.

Next we define $\P(n,m)$ to be 
\[
\{\,(p_1,\ldots,p_n,q_1,\ldots,q_m)\,|\,p_i>0,\,q_i\geq 0, \sum p_i +\sum
q_i=1\,\}
\]

Next we define the spaces $\C(n)$.  We define $\C(0)$ to be a point.
For $n>0$ the space $\C(n)$ that we are seeking to define is a quotient 
\[
\left(\bigcup_{m\geq 0} \P(n,m)\times \F(n,m)\right)/\sim
\] 
and our next task is to define the equivalence relation $\sim$. 
If $f\in\I_n^m$ and $1\leq k\leq m$ we write $f_k$ for the formula obtained 
from $f$ by ``pruning'' the $k$-th copy of $\varepsilon$ in $f$ (counting 
from left to right).  What this means is that if the copy of $\varepsilon$ 
has valence $j$ with $j>0$ we replace that $\varepsilon$ by $\bar{e}_j$ 
(see Remark \ref{hi}).  If the copy of $\varepsilon$ has no entries then we 
remove it, together with the typographical symbols immediately adjacent to it 
(these can be commas, parentheses, or $\smallsmile$'s). Now let $f\in\I_n^m$, 
and let
\[
\x=(\p,\q),(\s,\t)
\]
be a point of $\P(n,m)\times\F_f$. If $q_k=0$ then $\x$ will be 
equivalent to a point of $\P(n,m-1)\times\F_{f_k}$.  There are three cases.  If
the $k$-th copy of $\varepsilon$ in $f$ has valence 0 and if this $\varepsilon$
is the $i$-th entry of the integer $k'$ then $\x$ is
equivalent to the point
\[
(\p,q_1,\ldots,\widehat{q_k},\ldots,q_m),
(\s_1,\ldots,s^i\s_{k'},\ldots,\s_n,\t_1,\ldots,\widehat{t_k},\ldots,\t_m)
\]
of $\P(n,m-1)\times\F_{f_k}$  (where the hats indicate that the coordinates
$q_k$ and $\t_k$ are deleted, and $s^i$ is the $i$-th degeneracy map of the
simplex $\Delta^{v_f(k')})$. If the $k$-th copy of $\varepsilon$ in $f$ has
valence 0 and this copy of $\varepsilon$ is the $i$-th entry of the $k'$-th
copy of $\varepsilon$ in $f$ then $\x$ is equivalent to
\[
(\p,q_1,\ldots,\widehat{q_k},\ldots,q_m),
(\s,t_1,\ldots,s^i\t_{k'},\ldots,\widehat{t_k},\ldots,\t_m).
\]
Otherwise $\x$ is equivalent to
\[
(\p,q_1,\ldots,\widehat{q_k},\ldots,q_m),
(\s,t_1,\ldots,\widehat{t_k},\ldots,\t_m).
\]
This completes the definition of the equivalence relation $\sim$ and of the 
space $\C(n)$.

\begin{remark}
\label{hiyourself}
{\rm
For later use we remark that $\C'(n)$ is a deformation retract of $\C(n)$ by
the homotopy $H_t$ with 
\[
H_t((\p,\q),(\s,\t))=\left(\frac{1-t\sum q_i}{\sum p_i}\p,t\q\right),(\s,\t)
\]
}
\end{remark}

\bigskip

It remains to define the operad structure of $\C$ and the action of $\C$ on
${\mbox{\rm Tot}}(X^\b)$ when $X$ is the cosimplicial space associated to an operad with
multiplication. 

First we define the degeneracy maps of $\C$.  If $f\in\I_{n,m}$ and $1\leq 
k\leq n$ we write $_k f$ for the formula obtained from $f$ by replacing the 
symbol $k$ by $\varepsilon$, and then replacing $k+1,\ldots,n$ respectively 
by $k,\ldots,n-1$; thus $_k f$ is an element of $\I_{n-1}^{m+1}$.  
The degeneracy map
\[
\circ_k:\C(n)\rightarrow\C(n-1)
\]
takes a point
\[
(\p,\q),(\s,\t)
\]
of $\P(n,m)\times\F_f$ to the point
\[
((p_1,\ldots,\widehat{p_k},\ldots,p_n),(q_1,\ldots,p_k,\ldots,q_n)),
((\s_1,\ldots,\widehat{\s_k},\ldots,\s_n),(\t_1,\ldots,\s_k,\ldots,\t_n)),
\]
of $\P(n-1,m+1)\times\F_{_kf}$ (the hats mean that the coordinates $p_k$ 
and $\s_k$ are deleted).

Notice in particular that every point of $\C(n)$ is an iterated degeneracy of
some point in $\C'(n)$.  This implies that the operad structure maps of $\C$
and its action on ${\mbox{\rm Tot}}(X^\b)$ are determined by the corresponding data for
$\C'$, and it is straightforward to verify that the conditions for an operad
structure and an operad action are satisfied.
We are now finished with the definition of $\C$ and the proof of
part (ii) of Theorem \ref{main}.

\newpage

\section{The chain operad of $\C$.}
\label{sec6}

The (normalized) singular chains 
functor, applied to a topological operad, gives a chain operad since the
shuffle map is strictly associative and commutative. In this 
section we show that the singular chain operad of our operad $\C$ is 
quasi-isomorphic to the chain operad $\H$ described in the introduction.

We begin with a general result, and for this we need some terminology.  We will
consider filtered chain complexes of abelian groups (all of our chain 
complexes will be non-negatively graded, and all of our filtrations will be 
increasing).  We say that a chain complex $X$ is {\it homologically filtered} 
if the homology of $F_i X/F_{i-1}X$ is concentrated in dimension $i$.  (The
motivating example is when $X$ is the singular chain complex of a CW-complex 
$A$ and $F_iX$ is the singular chain complex of the $i$-skeleton of $A$.)  If 
$X$ is a homologically filtered chain complex we define the {\it 
condensation} of $X$, denoted $\overline{X}$, to be the chain complex whose 
$i$-th group is $H_i(F_iX/F_{i+1}X)$ and whose $i$-th boundary operator is 
the boundary operator of the triple $(F_{i}X,F_{i-1}X,F_{i-2} X)$. (In the 
motivating example, the condensation of $X$ is the usual cellular chain 
complex of $A$).  Note for later use that the homology of $\overline{X}$ is
naturally isomorphic to the homology of $X$ (this follows from the usual
spectral sequence argument).

By a {\it filtered chain operad} we mean a chain operad
$\D$ together with an increasing filtration $F_i$ of each $\D(n)$ such that 
for each $n$, $j$ and $k$ the operation
\[
\circ_k:\D(n)\otimes\D(j)\rightarrow \D(n+j-1)
\]
takes $F_i\D(n)\times F_{i'}\D(j)$ to $\F_{i+i'}\D(n+j-1)$.  
We say that $\D$ is {\it homologically filtered} if the filtration of each
$\D(n)$ is homological.  In this case we can define the {\it condensation} of
$\D$ to be the chain operad $\overline{\D}$ whose $n$-th object is the 
condensation $\overline{\D(n)}$ and whose operad structure maps are 
determined by those of $\D$.

\begin{theorem}
\label{jeff}
Let $\D$ be a homologically filtered chain operad.  Then $\D$ is
quasi-isomorphic in the category of chain operads to its condensation
$\overline{\D}$.
\end{theorem}

Before proving this we give some applications to our situation.  
Recall the chain operad $\H$ defined in the introduction.
Let $\H'$ be the suboperad of $\H$ defined by
$\H'(n)=\H(n)$ for $n>0$ and
$\H'(0)=0$.
We can filter $\C'$ by declaring $\P(n)\times\F_f$ to be in filtration $\sum
v_f(i)$.  It is easy to check that this induces a homological filtration of the
singular chains operad of $\C'$, and that the condensation of this filtered
chain operad is isomorphic to $\H'$. Now Theorem \ref{jeff} implies 

\begin{corollary}
The singular chain operad of $\C'$ is quasi-isomorphic to $\H'$. 
\end{corollary}

Similarly, we can filter $\C$ by declaring $\P(n,m)\times\F_f$ to be in 
filtration $\sum v_f(i)+\sum v_f(\varepsilon)$, and  this induces a homological 
filtration of the singular chain operad of $\C$.  Now the 
deformation retraction defined in Remark \ref{hiyourself} is filtration
preserving, and thus the inclusion $\C'(n)\hookrightarrow\C(n)$ induces an 
isomorphism (not just a quasi-isomorphism) of condensations for all $n>0$.  
Moreover, the definition of the degeneracy maps for $\C$ given in Section 
\ref{sec5} shows that the degeneracy maps in the condensation of $\C$ are 
zero (because if we begin with a cell in $\C(n)$, apply a degeneracy, and 
then apply the retraction of Remark \ref{hiyourself} we will end up in a 
lower filtration than that of the original cell). Together, these facts imply 
that the condensation of $\C$ is isomorphic to $\H$, and Theorem \ref{jeff} 
implies

\begin{corollary}
\label{wed2}
The singular chain operad of $\C$ is quasi-isomorphic to $\H$.
\end{corollary}

We now turn to the proof of Theorem \ref{jeff}.  

Let $Op$ denote the category
of non-negatively graded chain operads, let ${\mathbb N}$ be the non-negative
integers, and let $Ch^{\mathbb N}$ be the category whose objects are sequences of
chain complexes, indexed by $\mathbb N$, and whose morphisms are sequences of
chain maps.  The forgetful functor
\[
Op\rightarrow Ch^{\mathbb N}
\]
has a left adjoint, the free functor, which we will denote by 
\[
\Phi:
Ch^{\mathbb N} \rightarrow Op
\]

We will write $S^n$ for the chain complex consisting of a copy of $\mathbb Z$ in
dimension $n$, and $B^n$ for the chain complex 
\[
\cdots 0 \rightarrow {\mathbb Z} \labarrow{\id} {\mathbb Z} \rightarrow  0\cdots
\]
where the copies of $\mathbb Z$ are in dimensions $n$ and $n-1$ (of course, the
analogy is with the $n$-sphere and the $n$-disk).  Given a sequence of sets
$T=\{T_m\}_{m\geq 0}$, we will write $S^n(T)$ for the sequence of chain 
complexes
\[
\{\, \mathop{\oplus}_{T_m} S^n \,\}_{m\geq 0}
\]
and $B^n(T)$ for the sequence of chain complexes
\[
\{\, \mathop{\oplus}_{T_m} B^n \,\}_{m\geq 0}
\]

Now suppose we are given a chain operad $\D$, a sequence of sets $T$, 
and a sequence of chain maps
\[
f=\{\,f_m:\mathop{\oplus}_{T_m} S^{n-1} \rightarrow \D(m)\,\}
\]
Passing to the adjoint gives a map of chain operads
\[
\tilde{f}:\Phi (S^{n-1}(T))\rightarrow \D
\]
and we can form a pushout in the category of chain operads:
\[
\begin{array}{ccc}
\Phi(S^{n-1}(T))&\labarrow{\Phi{\iota}} & \Phi(B^n(T)) \\
\llap{$\tilde{f}$}\downarrow && \downarrow \\
\D &\rightarrow& \D'
\end{array}
\]
where $\iota$ is the inclusion $S^{n-1}(T)\hookrightarrow B^n(T)$.  
We then say that $\D'$ is obtained from $\D$ by {\it attaching $n$-cells}.
If we have a sequence of operads $\D_n$ such that $\D_0=S^0(T)$
for some $T$ and each $\D_n$ is
obtained from $\D_{n-1}$ by attaching $n$-cells, we say that the colimit
$\displaystyle\mathop{\mbox{colim}}_n \D_n$ is a {\it CW chain operad}.

\begin{lemma} \label{jeffone}
If $\D$ is any chain operad, there is a CW chain operad $\E$ and a
quasi-isomorphism $\E\rightarrow \D$.
\end{lemma}

The proof of the lemma is precisely analogous to that of the corresponding fact
in the category of spaces, so we omit it.

We need to define what it means for two morphisms of chain operads to be
chain-homotopic ``through operad maps.''  First note that if $C$ is
a differential graded coalgebra and $\D$ is any chain operad we can define 
a new chain operad $\Hom(C,\D)$ by letting the $n$-th object be 
$\Hom(C,\D(n))$ and letting the operad structure be determined by that of 
$\D$.  Now let $I$ be the simplicial chain complex of the standard 1-simplex.
Let $i_0$ (respectively $i_1$) be the maps from $S^0$ to $I$ corresponding to
the endpoints of the 1-simplex.  Then two chain-operad morphisms 
$f_1,f_2:\D\rightarrow\D'$ are {\it operad-chain-homotopic} if there is a 
chain-operad morphism $H:\D\rightarrow\Hom(I,\D')$ with $i_1^*\circ H=f_1$ 
and $i_2^*\circ H=f_2$.

Next let us observe that if $C$ is a chain complex, then we can give $C$
the filtration $F_i$ where $F_i(C)$ is equal to $C$ in dimensions $\leq i$
and is 0 in dimensions $>i$.  We call this the {\it filtration by degrees}.  
This is a homological filtration, and the associated condensation is 
isomorphic to the $C$ that we started with.

\begin{proposition} \label{jefftoo}
Let $\E$ be a CW chain operad, $\D$ a homologically filtered chain operad, 
and $f:\E\rightarrow \D$ a morphism of chain operads.   Give $\E$ the 
filtration by degrees. Then $f$ is operad-chain-homotopic to a morphism $f'$ 
of filtered chain operads.
\end{proposition}

Before proving Proposition \ref{jefftoo} we use it to prove Theorem
\ref{jeff}.  So let $\D$ be a chain operad with a homological filtration. By
Lemma \ref{jeffone} there is a CW chain operad $\E$ and a quasi-isomorphism of
chain operads $f:\E\rightarrow \D$.  By Proposition \ref{jefftoo}, $f$ is
operad-chain-homotopic to a morphism $f'$ of filtered chain operads (where $\E$
is given the filtration by degrees.)  Now $f'$ induces an operad morphism of 
condensations 
\[
\overline{f'}:\overline{\E}\rightarrow \overline{\D},
\]
and $\overline{f'}$ is a quasi-isomorphism (because $\overline{f'}$ has the 
same effect in homology as $f'$, and $f'$ has the same effect in homology as
$f$).  But $\overline{\E}$ is isomorphic to $\E$ (since $\E$ was given the
filtration by degrees), so we conclude that $\E$ is quasi-isomorphic to
$\overline{\D}$.  Since we also know that $\E$ is quasi-isomorphic to
$\D$, we have finally that $\D$ is quasi-isomorphic to $\overline{\D}$, as
required.
\qed

It remains to prove Proposition \ref{jefftoo}.  First we need two lemmas.

\begin{lemma}\label{sat1}
Let $\D$ be a filtered operad, let $T$ be a sequence of sets, and for each 
$m$ let $f_m:S^n(T_m)\rightarrow \D(m)$ be a chain map which lands in 
filtration $n$.  Let $f$ denote the sequence $\{f_m\}$. Then the chain-operad 
map
\[
\tilde{f}:\Phi(S^n(T)) \rightarrow \D
\]
is filtration preserving when $\Phi(S^n(T)$ is given the filtration by
degrees.
\end{lemma}

\noindent
{\bf Proof.}
This is immediate from the explicit description of $\Phi$, which may be
found in \cite{GJ}.
\qed

\begin{lemma}\label{sat2}
Let 
\[
\begin{array}{ccc}
\D_0&\rightarrow&\D_1 \\
\downarrow&&\downarrow \\
\D_2&\rightarrow&\D
\end{array}
\]
be a pushout diagram in the category of chain operads.  Suppose that $\D'$ is a
filtered chain operad and that
$f:\D\rightarrow \D'$ is a morphism of chain operads.  Give $\D_1$, $\D_2$ and
$\D$ the filtration by degrees and suppose that $f$ is filtration-preserving
when restricted to $\D_1$ and $\D_2$.  Then $f$ is filtration-preserving. 
\end{lemma}

\noindent
{\bf Proof.}
Suppose first that the pushout diagram is obtained by applying $\Phi$ to a
pushout diagram in the category $Ch^{\mathbb N}$.  In this case the lemma follows
from Lemma \ref{sat1} by passage to adjoints.  But the general case can be
reduced to this case by using the fact that the pushout of 
\[
\begin{array}{ccc}
\Phi\D_0&\rightarrow&\Phi\D_1 \\
\downarrow&&\\
\Phi\D_2&&
\end{array}
\]
surjects to $\D$.
\qed

\medskip

\noindent
{\bf Proof of Proposition \ref{jefftoo}.}
Let $\E_n$ be a sequence of chain operads, with $\E_n$ obtained from $\E_{n-1}$
by attaching $n$-cells, and let $\E=\displaystyle\mathop{\mbox{colim}}_n 
\E_n$. Suppose inductively that we are given a homotopy 
$H_{n-1}:\E_{n-1}\rightarrow\Hom(I,\D)$ such that $i_0^*\circ H_{n-1}=f$
and $i_1^*\circ H_{n-1}$ is filtration-preserving.  It suffices to show that
$H_{n-1}$ extends to a homotopy $H_n:\E_{n}\rightarrow\Hom(I,\D)$ such that 
$i_0^*\circ H_{n}=f$ and $i_1^*\circ H_{n}$ is filtration-preserving.  
Now $\E_n$ is defined by a pushout diagram
\[
\begin{array}{ccc}
\Phi(S^{n-1}(T))&\labarrow{\Phi{\iota}} & \Phi(B^n(T)) \\
\llap{$\tilde{g}$}\downarrow && \downarrow\rlap{$\gamma$} \\
\E_{n-1} &\rightarrow& \E_n
\end{array}
\]
and it suffices (by Lemma \ref{sat2}) to show that $H_{n-1}\circ\tilde{g}$ 
extends to a homotopy
\[
H':\Phi(B^n(T))\rightarrow \Hom(I,\D)
\]
such that $i_0^*\circ H'=f\circ\gamma$ and $i_1^*\circ H'$ is
filtration-preserving. 

Passing to adjoints, we need to find a suitable chain homotopy 
\[
s:B^n(T)\rightarrow\D,
\]
(where $s$ raises degrees by 1) and we can treat each of the summands $B^n$ 
in $B^n(T)$ separately, so we are in the following situation.  Choose 
generators $a$ and $b$ for the two copies of $\mathbb Z$ in $B^n$ such that 
$\partial b=a$.  The homotopy $H_{n-1}$ determines $s(a)$ and we have
\[
\partial s(a)=f(a)-x,
\]
where $x$ is in filtration $n-1$.  We need to define $s(b)$ in such a way that
\[
\partial(s(b))+s(a)=f(b)-y,
\]
where $y$ is in filtration $n$ and $\partial y=x$.  We know that 
\[
\partial f(b)=f(a),
\]
and this implies that
\[
\partial (f(b)-s(a))=x.
\]
The fact that $\D$ is homologically filtered now implies that there is an
element $z$ in filtration $n$ with 
\[
\partial z= x.
\]
Since $\partial z=\partial(f(b)-s(a))=x$, we see that $f(b)-s(a)-z$ is a cycle,
so (using again the fact that $\D$ is homologically filtered) there is a cycle
$z'$ in filtration $n$ which is homologous to $f(b)-s(a)-z$; that is, there is
an element $w$ in dimension $n+1$ with
\[
\partial w=(f(b)-s(a)-z)-z'
\]
We can now define $s(b)=w$ and let $y=z+z'$.
\qed

\newpage

\section{The equivalence between $\C$ and the little 2-cubes operad.}
\label{sec7}

In this section and the next, we will prove

\begin{theorem} \label{mon2}
$\C$ and $\C_2$ are weakly equivalent operads.
\end{theorem}

For the proof we use a method of
Fiedorowicz \cite{Fied}.  
Let $\C_1$ denote the operad of little
intervals and let $\tau_i\in\Sigma_n$ denote the permutation that 
transposes $i$ and $i+1$.

\begin{proposition}
\label{mon1}
Suppose we are given an operad $\D$ with
\begin{description}
\item[{\rm (a)}]
a morphism of (non-$\Sigma$) operads $I: \C_1\rightarrow\D$, 
\item[{\rm (b)}] for each $n$, a point $c_n\in \C_1(n)$ and for each $i$ a 
path 
$\alpha_i$ from $I(c_n)$ to $\tau_i I(C_n)$. 
\end{description}
Suppose moreover that 
\begin{description}
\item[{\rm (c)}] the universal cover of $\D(n)$ is contractible for each $n$,
and
\item[{\rm (d)}] for each $n$ and $i$ the paths 
\[
\tau_i\tau_{i+1}(\alpha_i)\cdot \tau_i(\alpha_{i+1})\cdot \alpha_i
\]
and
\[
\tau_{i+1}\tau_{i}(\alpha_{i+1})\cdot \tau_{i+1}(\alpha_{i})\cdot \alpha_{i+1}
\]
are path homotopic (where $\cdot$ denotes concatenation of paths).
\end{description}
Then $\D$ is weakly equivalent as an operad to $\C_2$ (in the sense that there
is a third operad which maps to each of them by a morphism which is a weak
equivalence on each space.)
\end{proposition}

\bigskip
\noindent
{\bf Proof of Proposition \ref{mon1}} (following Fiedorowicz).
For each $n$ let $\tilde{\D}(n)$ be the universal cover of the space $\D(n)$,
and let $\pi:\tilde{\D}(n)\rightarrow\D(n)$ be the projection.  We assume for
simplicity that $I:\C_1(n)\rightarrow \D(n)$ is an inclusion.  Since $\C_1(n)$
is contractible,
$\pi^{-1}\C_1(n)$ is a disjoint union of copies of $\C_1(n)$; we arbitrarily
choose one such copy for each $n$ and declare that the basepoint of
$\tilde{\D}(n)$ shall lie in it.  Now that we have chosen basepoints,
$\tilde{\D}$ has an induced structure of non-$\Sigma$ operad.  Next observe
that hypothesis (b) determines a map
$\tilde\tau_i:\tilde\D(n)\rightarrow\tilde\D(n)$ for each $n$ and $i$, and 
hypothesis (d) implies that the braid relation
\[
\tilde\tau_i\tilde\tau_{i+1}\tilde\tau_i=\tilde\tau_{i+1}\tilde\tau_i
\tilde\tau_{i+1}
\]
is satisfied, so each $\tilde\D(n)$ has an action of the braid group $B_n$.
This makes $\tilde\D$ into a braided operad as defined by Fiedorowicz (this
just means that the symmetric groups are replaced by the braid groups
everywhere in the usual definition of operad).  Next let $\tilde\C_2$ be the
corresponding braided operad constructed from $\C_2$ (see \cite[Example 
3.1]{Fied}).  Since $\C_2(n)$ is a
$K(B_n,1)$ for each $n$, the spaces $\tilde\C_2(n)$ are contractible for each
$n$, so by hypothesis (c) the projections
\[
\tilde\D\times\tilde\C_2 \rightarrow \tilde\D
\]
and
\[
\tilde\D\times\tilde\C_2 \rightarrow \tilde\C_2
\]
are equivalences of braided operads.  If we now mod out by the action of the
pure braid groups (that is, the kernels of the projections
$B_n\rightarrow\Sigma_n)$, we get the assertion of the proposition.
\qed.

Now we turn to the proof of Theorem \ref{mon2}.
We need to show that our operad $\C$ satisfies the hypotheses of Proposition
\ref{mon1}. For (a), let $c\in\C_1(n)$.  Let $(p_1,\ldots,p_n$ be the lengths
of the little intervals in $c$ and let $(q_1,\ldots,q_{n+1})$ be the lengths of
the gaps, and
define $I(c)=(\p,\q,e\smallsmile 1\smallsmile e \smallsmile \ldots 
\smallsmile n \smallsmile e)$. 
For (b), we define $c_n$ to consist of $n$ intervals each of length
$\frac{1}{n}$,
and we define $\alpha_i$ to be the 1-simplex
\[
\{ c_n\}\times \F_{1\smallsmile\ldots\smallsmile (i-1)\smallsmile i(i+1) 
\smallsmile (i+2)\smallsmile \ldots \smallsmile n}
\]
We defer the verification of (c) to the next section. For (d), assume first
that $n=3$ and $i=1$.  In this case the desired path homotopy is given by the
following picture: 

\bigskip
\bigskip
\bigskip
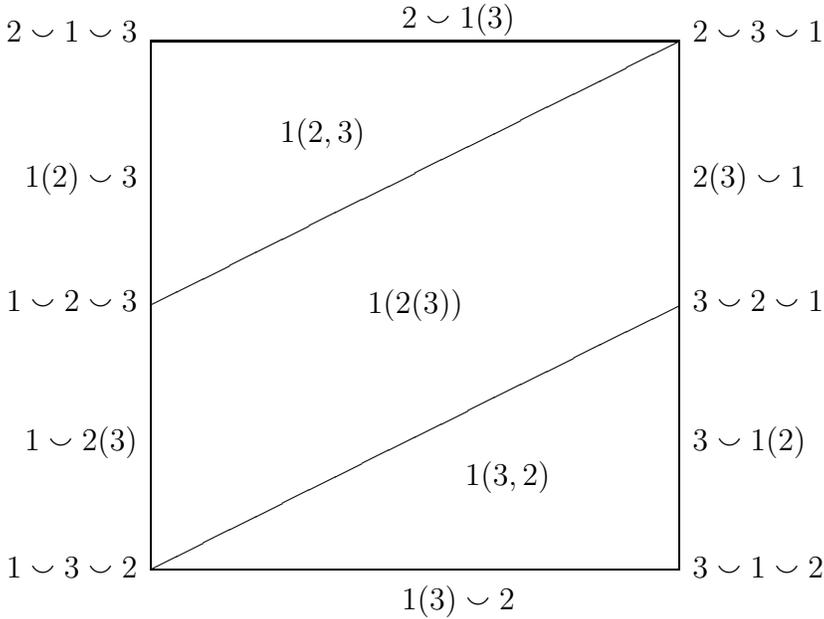
\begin{picture}(200,200)

\put(0,0){\line(1,0){200}}
\put(0,0){\line(0,1){200}}
\put(0,200){\line(1,0){200}}
\put(200,0){\line(0,1){200}}

\put(0,0){\line(2,1){200}}
\put(0,100){\line(2,1){200}}

\put(0,95){\makebox(200,10){$1(2(3))$}}

\put(0,160){\makebox(130,10){$1(2,3)$}}
\put(70,30){\makebox(130,10){$1(3,2)$}}

\put(-5,45){\llap{$1\smallsmile 2(3)$}}
\put(205,45){\rlap{$3\smallsmile 1(2)$}}

\put(-5,145){\llap{$1(2)\smallsmile 3$}}
\put(205,145){\rlap{$2(3)\smallsmile 1$}}

\put(-5,98){\llap{$1\smallsmile 2\smallsmile 3$}}
\put(205,98){\rlap{$3\smallsmile 2\smallsmile 1$}}

\put(95,-15){\rlap{$1(3)\smallsmile 2$}}
\put(95,205){\rlap{$2  \smallsmile 1(3)$}}

\put(-5,-3){\llap{$1\smallsmile 3\smallsmile 2$}}
\put(-5,200){\llap{$2\smallsmile 1\smallsmile 3$}}
\put(205,-3){\rlap{$3\smallsmile 1\smallsmile 2$}}
\put(205,200){\rlap{$2\smallsmile 3\smallsmile 1$}}
\end{picture}

\bigskip
\bigskip
\bigskip

Here the path
\[
\tau_1\tau_{2}(\alpha_1)\cdot \tau_1(\alpha_{2})\cdot \alpha_1
\]
begins at the vertex $1\smallsmile 2\smallsmile 3$ and goes clockwise to 
$3\smallsmile 2\smallsmile 1$ while the path
\[
\tau_{2}\tau_{1}(\alpha_{2})\cdot \tau_{2}(\alpha_{1})\cdot \alpha_{2}
\]
begins at the vertex $1\smallsmile 2\smallsmile 3$ and goes counterclockwise to 
$3\smallsmile 2\smallsmile 1$. 

The verification of (d) for general $n$ and $i$ is similar: the required
homotopy is again the boundary of a figure consisting of a square and two
triangles, but the triangle labeled  $1(2,3)$ in the figure above is
replaced by 
\[
1\smallsmile\ldots \smallsmile (i-1)\smallsmile
i(i+1,i+2)\smallsmile (i+3) \smallsmile \ldots \smallsmile n
\]
and so on.

\newpage

\section{The homotopy type of the spaces $\C(n)$.}
\label{sec8}

In this section we prove 

\begin{proposition} \label{wed4} 
For each $n\geq 0$ the space $\C(n)$ is homotopy
equivalent to the space $\C_2(n)$.  
\end{proposition}

On passage to universal covers this will
imply that hypothesis (c) of Proposition \ref{mon1} is satisfied for the operad
$\C$, and this will complete the proof of Theorem \ref{mon2}.

Since $\C(0)$ is a point we may assume $n>0$.

The homotopy in Remark \ref{hiyourself} shows that $\C(n)$ is homotopy
equivalent to $\F(n)$, and it is well-known that $\C_2(n)$ is homotopy
equivalent to the configuration space $F(n)$ of $n$ ordered points in 
${\mathbb R}^2$, so what we really need to show is that $\F(n)$ is weakly 
equivalent to 
$F(n)$.

The basic idea is to show that $\F(n)$ and $F(n)$ are both weakly
equivalent to the nerve of a certain category $\T_n$.  In order to
motivate the definition of $\T_n$, let us observe that if $f$ is a
formula of type $n$ then $f$ induces a total order on the set
$\{1,\ldots,n\}$ (namely the order in which these symbols first appear
in $f$ as one reads from left to right) and also a partial order
(generated by the relation: $i< j$ if $j$ is an entry of $i$).  These
two orders are consistent in the sense that
\begin{equation}
\label{tu1}
\mbox{$p\subseteq t$, and if $i< j< k$ in $t$ and $i< k$ in $p$ then $i$ must 
be $< j$ in $p$.}
\end{equation}
(Here $p\subseteq t$ means that $p$ is contained in $t$ when both are
considered as sets of ordered pairs, i.e., if $i< j$ in $p$ then $i< j$ in $t$.)

Accordingly, we define $\T_n$ to be the following partially ordered set.
The objects $\T_n$ are pairs $(t,p)$, where $t$ is a total order of the set 
$\{1,\ldots,n\}$ and $p$ is a partial order of this set, subject to the
consistency condition (\ref{tu1}).  In order to define the partial ordering of
$\T_n$, let us first define $t^{op}$ to be the ordering which is the
reverse of $t$ (that is, $i< j$ in $t^{op}$ if and only if $j< i$ in $t$).
Then there is a morphism in $\T_n$ from $(t_1,p_1)$
to $(t_2,p_2)$ if both of the conditions $p_1\subseteq p_2$ and $t_2\cap 
t_1^{op}\subseteq p_2$ are satisfied.
(We don't know a good way to motivate this partial order on $\T_n$, except that
it is what is needed to make the proof work. But note that if $f< f'$ in the
partial order of subsection \ref{subseca2}, then the pair $(t,p)$ associated to 
$f$ is $<$ the pair associated to $f'$ in the partial order of $\T_n$.) 

The main step in proving Proposition \ref{wed4} is to show

\begin{proposition} \label{tu2}
\begin{description}
\item[{\rm (a)}]
There is a functor $I_n$ from $\T_n$ to spaces such that
\begin{description}
\item[{\rm (i)}]
the spaces $I_n(t,p)$ are
an open cover of $F(n)$,
\item[{\rm (ii)}]
$F(n)$ is the colimit of $I_n$, and
\item[{\rm (iii)}]
each $I_n(t,p)$ is contractible.
\end{description}
\item[{\rm (b)}]
There is a functor $I'_n$ from $\T_n$ to spaces such that
\begin{description}
\item[{\rm (i)}]
the spaces $I'_n(t,p)$ are
subcomplexes of $\F(n)$,
\item[{\rm (ii)}]
$\F(n)$ is the colimit of $I'_n$, and
\item[{\rm (iii)}]
each $I'_n(t,p)$ is contractible.
\end{description}
\end{description}
\end{proposition}

\noindent
{\bf Proof of Proposition \ref{wed4}.}
Parts (a)(i) and (a)(ii), together with 
the proof of \cite[Proposition 4.1]{Segal}, imply that the natural map 
\[
{\rm hocolim} \, I_n\rightarrow {\rm colim}\, I_n
\]
is a weak equivalence.  Part (a)(iii), together with the homotopy invariance of
homotopy colimits (\cite[XII.4.2]{BK}), implies that hocolim$\, I_n$  is weakly
equivalent to hocolim$\,*$ (where $*$ is the constant functor that takes every
object of $\T_n$ to a point.)

Similarly, parts
(b)(i) and (b)(ii), together with 
\cite[Proposition 6.9]{Balt}, imply that the natural map 
\[
{\rm hocolim} \, I'_n\rightarrow {\rm colim}\, I'_n
\]
is a weak equivalence, and part (b)(iii)
implies that hocolim$\, I'_n$  is weakly
equivalent to hocolim$\,*$.
\qed

\medskip

We proceed to the proof of Proposition \ref{tu2}(a).  Let $\pi_1$ and $\pi_2$
denote the projections of ${\mathbb R}^2$ on its two factors.
Let $I_n(t,p)$ be the
open subspace of $F(n)$ consisting of all configurations $\x=(x_1,\ldots,x_n)$
such that 
\begin{description}
\item \qquad
if $i< j$ in $t$ and $\pi_1(x_i)\geq\pi_1(x_j)$ then $i< j$ in $p$
\item \qquad
if $i< j$ in $t$ and $\pi_1(x_i)=\pi_1(x_j)$ then $i< j$ in $p$ and
$\pi_2(x_i)<\pi_2(x_j)$.
\end{description}

Part (a)(i) is evident.  For part (a)(ii), since the sets 
$I_n(t,p)$ cover $F(n)$ it suffices to show that the map 
colim$\,I_n \rightarrow F(n)$ is 1-1, and for this it suffices to show that, 
if $\x$ is in $I_n(p,t)\cap I_n(p',t')$, there is a pair $(p'',t'')$ which is 
$\leq$ both $(p,t)$ and $(p',t')$ and is such that $\x\in I_n(p'',t'')$.  And 
this in turn follows from the stronger fact that for each configuration $\x$ 
there is a unique minimal $(t,p)$ with $\x\in I_n(t,p)$: define $t$ by  
$i< j$ in $t$ if and only if $x_i<x_j$ in lexicographic order, and define 
$p$ by $i< j$ in $p$ if and only if $\pi_1(x_i)=\pi_1(x_j)$ and 
$\pi_2(x_i)<\pi_2(x_j)$.

The proof of part (a)(iii) is by induction on $n$. The case $n=1$ is obvious,
so assume $n>1$, and let $(t,p)\in\T_n$.  To simplify the notation we treat 
the case where $t$ is the standard total order $1<2<\ldots<n$; the general 
case is precisely similar.  Let $t'$ be the standard total order of
$\{1,\ldots,n-1\}$, and let $p'$ be the restriction of $p$ to 
$\{1,\ldots,n-1\}$. It suffices by induction to show that $I_n(t,p)$ is 
homotopy equivalent to $I_{n-1}(t',p')$.  Let 
\[
\alpha:I_n(t,p)\rightarrow I_{n-1}(t',p')
\]
be the projection map that takes $(x_1,\ldots,x_n)$ to $(x_1,\ldots,x_{n-1})$.
Let 
\[
\beta:I_n(t,p)\rightarrow I_{n-1}(t',p')
\]
take $(x_1,\ldots,x_{n-1})$ to
$(x_1,\ldots,x_{n-1},\gamma(x_1,\ldots,x_{n-1}))$,
where 
\[
\gamma(x_1,\ldots,x_{n-1})=(1+\max\{\pi_1(x_i)\},1+\max\{\pi_2(x_i)\})\in
{\mathbb R}^2
\]
Then $\alpha\circ\beta$ is the identity map, and it suffices to show that
$\beta\circ\alpha$ is homotopic to the identity map of $I_n(t,p)$.  
Let $H$ be the homotopy which leaves $x_1$,\ldots,$x_{n-1}$ fixed and moves
$x_n$ vertically until it has the desired 
second coordinate, and then horizontally until it has the desired first 
coordinate: that is,
\[
H_t(\x)=\left\{
\begin{array}{l}
(x_1,\ldots,x_{n-1}, (1-2t)x_n+2t(\pi_1(x_n),\pi_2(\gamma)))
\quad\mbox{if $t\leq 1/2$} \\
(x_1,\ldots,x_{n-1}, (2-2t)(\pi_1(x_n),\pi_2(\gamma)) +(2t-1)\gamma)
\quad\mbox{if $1/2\leq t$} \\
\end{array}
\right.
\]
where we have written $\gamma$ for $\gamma(x_1,\ldots,x_{n-1})$.
(To see that this homotopy stays inside $I_n(t,p)$ note that there cannot be 
a point in $\x$ directly above $x_n$, and if $x_i$ is to the right of $x_n$ we 
must have $i< n$ in $p$.) This concludes the proof of part (a) of 
Proposition \ref{tu2}.

For part (b), let $I'_n(t,p)$ be the following subcomplex of $\F(n)$:
\[
I'_n(t,p)=\bigcup_{\{f\,|\,(t_f,p_f)\leq (t,p)\}} \F_f
\]
where $(t_f,p_f)$ is the pair determined by the formula $f$ and $\F_f$ is as 
in subsection \ref{subseca1}. 

Part (b)(i) is evident.  For part (b)(ii), since the sets 
$I'_n(t,p)$ cover $\F(n)$ it suffices to show that the map 
colim$\,I'_n \rightarrow \F(n)$ is 1-1, and for this it suffices to show that, 
if $\F_f$ is contained in $I'_n(p,t)\cap I'_n(p',t')$, there is a pair 
$(p'',t'')$ which is $\leq$ both $(p,t)$ and $(p',t')$ and is such that 
$\F_f\subseteq I'_n(p'',t'')$.  And 
this in turn follows from the stronger fact that for each formula $f$
the pair $(t_f,p_f)$ is
the unique minimal pair whose image under $I_n$ contains $\F_f$.

We now turn to the proof of (a)(iii), so let us fix a pair $(t,p)$.  
If $i\in\{1,\ldots,n\}$, we write $h(i)$ (the height of $i$ with respect to 
$p$) for the length of the longest chain $j_1<\ldots<i$ in $p$, and we write
$h(p)$ (the height of $p$) for the length of the longest chain in $p$. 
We write $w(p)$ (the width of $p$) for the number of elements 
$i\in\{1,\ldots,n\}$ with $h(i)=h(p)$.  The proof of (a)(iii) will be by double
induction on $h(p)$ and $w(p)$.

Let us fix an element $i$ with $h(i)=h(p)-1$.
There are two cases: 
\begin{description}
\item[Case 1] There is only one $j$ with $i<j$ in $p$.
\item[Case 2] There is more than one $j$ with $i<j$ in $p$.
\end{description}

We begin with Case 1.  In this case, condition (\ref{tu1}) implies that $j$ is
the immediate successor of $i$ in the total order $t$.

We define three subcomplexes of $I_n(p,t)$.  It will be convenient to write 
$f\leq (t,p)$ to mean $(t_f,p_f)\leq (t,p)$.

\begin{description}
\item $A=\bigcup_{f\in S_1} \F_f$, where $S_1=\{\,f\leq (t,p)\,|\, 
\mbox{$f$ contains one of the strings $i(j)$, $i\smallsmile j$ or
$j\smallsmile i$}\}$.
\item $B=\bigcup_{f\in S_2} \F_f$, where $S_2=\{\,f\leq (t,p)\,|\, 
\mbox{$i$ is to the left of $j$ in $f$ but $j$ is not an entry of $i$ in 
$f$}\}$
\item $C=\bigcup_{f\in S_3} \F_f$, where $S_3=\{\,f\leq (t,p)\,|\, 
\mbox{$i$ is to the right of $j$ in $f$}\}$.
\end{description}

Clearly $I_n(t,p)=A\cup B\cup C$.  We will show that $I_n(t,p)$ is contractible
by showing that $A$, $B$, $C$, $A\cap B$, and $A\cap C$ are contractible and
that $B\cap C$ is empty.

It is easy to see that $B=I_n(t,p-(i,j))$ (where $p-(i,j)$ means the
partial order which is the same as $p$ except that $i$ is no longer $<j$; this
is indeed a partial order since $j$ is the immediate successor of $i$ in $p$)
so $B$ is contractible by the inductive hypothesis.  Similarly,
$C=I_n(\bar{t},p-(i,j))$ (where $\bar{t}$ is the same as $t$ except that
$i$ and $j$ are switched; this is indeed a total order since $j$ is the
immediate successor of $i$ in $t$) so $C$ is contractible by the inductive
hypothesis.  It is also clear that $B\cap C$ is empty.

Next we claim that $A$ is contractible.  Let $_j t$ (respectively $_jp$) be 
the restriction of $t$ (respectively, $p$) to $\{1,\ldots,j-1,j+1,\ldots,n\}$.
If $i$ has valence $0$ in $g$ let us write $g*_i h$ to mean ``plug $h$ in for
$i$'' (this is a little different from the use of this symbol in subsection
\ref{subsecc}).
Then $S_1=\{ g*_i h\,|\, \mbox{$g\leq (_j t,_j p)$ and $h\leq i(j)$}\}$. This
implies that $A=I_{n-1}(_j t,_j p)\times \Delta^1$, and so $A$ is contractible
by the inductive hypothesis.

Similarly, $S_1\cap S_2=
\{ g*_i h\,|\, \mbox{$g\leq (_j t,_j p)$ and $h=i\smallsmile j$}\}$.  This
implies that
$A\cap B=I_{n-1}(_jt,_jp)\times \Delta^0$ which is contractible by the
inductive hypothesis. 
Also, $S_1\cap S_3=
\{ g*_i h\,|\, \mbox{$g\leq (_j t,_j p)$ and $h=j\smallsmile i$}\}$.  This
implies that
$A\cap B=I_{n-1}(_jt,_jp)\times \Delta^0$ which is contractible by the
inductive hypothesis.  This concludes the proof of Case 1.

For Case 2, let $j_1,\ldots,j_m$ be the successors of $i$ in $p$. We may assume
that $j_1<j_2<\ldots<j_m$ in $t$.  Condition (\ref{tu1}) then implies that 
$i,j_1,\ldots,j_m$ is a consecutive sequence in $t$ (that is, $j_1$ is the
immediate successor of $i$, etc.).

Let $t'$ (respectively $p'$) denote the restriction of $t$ (respectively $p$)
to $\{1,\ldots,n\}-\{j_1,\ldots,j_m\}$.

Now define subcomplexes $A$, $B$, and $C$ of $I_n(t,p)$ by

\begin{description}
\item $A=\bigcup_{f\in S_1} \F_f$, where $S_1=\{ 
g*_i h \,|\,
\mbox{$g\leq (t',p')$ and $h\leq i(j_1,\ldots,j_m)$
} \}$. 
\item $B=\bigcup_{f\in S_2} \F_f$, where $S_2=\{\,f\,|\,
\mbox{$j_1$ is to the left of $i$ in $f$}
\}$.
\item $C=\bigcup_{f\in S_3} \F_f$, where $S_3=\{ \,f\,|\,
\mbox{$j_m$ is to the right of $i$ in $f$ but is not an entry of $i$}
\}$.
\end{description}

Then $A\cup B\cup C=I_n(p,t)$, because if $f$ is not in $S_2$ or $S_3$ then 
$j_1$ will be to the right of $i$, $j_m$ will be an entry of $i$, and this will
imply that $f$ contains the string $i(j_1,\ldots,j_m)$ so that $f\in S_1$.

We will show that $I_n(t,p)$ is contractible by showing that $A$, $B$, $C$, 
$A\cap B$, $A\cap C$, $B\cap C$ and $A\cap B\cap C$ are all contractible.

First of all, $B$ is contractible by induction because it is equal to
$I_n(\bar{t},p-(i,j_1))$ (where $\bar{t}$ is the same as $t$
but with $i$ and $j_1$ switched).  Similarly, $C$ is contractible because it is
equal to $I_n(t,p-(i,j_m))$, and $B\cap C$ is contractible because it is equal
to $I_n(\bar{t},p-(i,j_1)-(i,j_s))$.

Next, $A$ is contractible because it is homeomorphic to
$I_{n-m}(t',p')\times\Delta^m$. Similarly, we have
\begin{itemize}
\item $S_1\cap S_2 =\{ 
g*_i h \,|\,
\mbox{$g\leq (t',p')$ and $h\leq j_1\smallsmile i(j_2,\ldots,j_m)$
} \}$, and so 
$A\cap B \approx I_{n-m}(t',p')\times \Delta^{m-1}$.
\item $S_1\cap S_3 =\{ 
g*_i h \,|\,
\mbox{$g\leq (t',p')$ and $h\leq i(j_1,\ldots,j_{m-1})\smallsmile j_m$
} \}$, and so 
$A\cap C \approx I_{n-m}(t',p')\times \Delta^{m-1}$.
\item $S_1\cap S_2 \cap S_3 =\{ 
g*_i h \,|\,
\mbox{$g\leq (t',p')$ and $h\leq j_1\smallsmile 
i(j_2,\ldots,j_{m-1})\smallsmile j_m$
} \}$, and so 
$A\cap B\cap C \approx I_{n-m}(t',p')\times \Delta^{m-2}$.
\end{itemize}

This concludes the proof of Case 2.\qed

\end{document}